\documentclass{dcds-b}
\usepackage{amsmath}
  \usepackage{paralist}
  \usepackage{graphics} %% add this and next lines if pictures should be in esp format
  \usepackage{epsfig} %For pictures: screened artwork should be set up with an 85 or 100 line screen
\usepackage{graphicx}  \usepackage{epstopdf}%This is to transfer .eps figure to .pdf figure; please compile your paper using PDFLeTex or PDFTeXify.
 \usepackage[colorlinks=true]{hyperref}
\hypersetup{urlcolor=blue, citecolor=red}

\def\currentvolume{20}
 \def\currentissue{8}
  \def\currentyear{2015}
   \def\currentmonth{October}
    \def\ppages{XX--XX}
     \def\DOI{10.3934/dcdsb.2015.20.xx}

\newtheorem{theorem}{Theorem}[section]

\newtheorem{assumption}{Assumption}

\theoremstyle{definition}
\newtheorem{definition}[theorem]{Definition}
\newtheorem{remark}{Remark}

\title[Classical Converse Lyapunov Theorems] %Use the shortened version of the full title
      {Classical Converse Theorems in Lyapunov's\\ Second Method}

\author[Christopher M. Kellett]{}


 \email{Chris.Kellett@newcastle.edu.au}
\thanks{C.~M.~Kellett is supported by Australian Research Council Future Fellowship FT1101000746.}

\usepackage{amssymb}

\usepackage{color}
\def\rev{\color{black}}

\usepackage[OT2,OT1]{fontenc} \newcommand\cyr{%
\renewcommand\rmdefault{wncyr}%
\renewcommand\sfdefault{wncyss}%
\renewcommand\encodingdefault{OT2}%
\normalfont
\selectfont}
\DeclareTextFontCommand{\textcyr}{\cyr}

\def\cprime{\char126}

\def\K{\mathcal{K}}
\def\Kinf{\K_{\infty}}
\def\KL{\mathcal{KL}}
\def\L{\mathcal{L}}

\def\R{\mathbb{R}}

\def\G{\mathcal{G}}

\def\S{\mathcal{S}}
\def\Z{\mathbb{Z}}
\def\O{\mathcal{O}}
\def\B{\mathcal{B}}
\def\A{\mathcal{A}}
\def\D{\mathcal{D}}

\begin{document}
\maketitle

% Enter the first author's name and address:
\centerline{\scshape Christopher M. Kellett }
\medskip
{\footnotesize
% please put the address of the first author
 \centerline{School of Electrical Engineering and Computer Science}
   \centerline{University of Newcastle}
   \centerline{ Callaghan, New South Wales, 2308, Australia}
} % Do not forget to end the {\footnotesize by the sign }

%\bigskip

% The name of the associate editor will be entered by an editorial staff
% "Communicated by the associate editor name" is not needed for special issue.
% \centerline{(Communicated by the associate editor name)}

%The abstract of your paper
\begin{abstract}
%This is the abstract of your paper and it should not exceed
%\textbf{200} words.
	Lyapunov's second or direct method is one of the most widely used techniques for
	investigating stability properties of dynamical systems.  This technique makes use of an
	auxiliary function, called a Lyapunov function, to ascertain stability properties for a specific
	system without the need to generate system solutions.  An important question is the
	converse or reversability of Lyapunov's second method; i.e., given a specific stability property
	does there exist an appropriate Lyapunov function?  We survey some of the available answers to this
	question.
\end{abstract}

\section{Introduction}
\label{sec:Intro}
Over the last 100 years, Lyapunov's second, or direct, method has arguably been the most widely used technique for
analyzing stability properties of various types of mathematically described dynamical systems, including differential
and difference equations, hybrid differential-difference equations, stochastic differential equations, and many others.
In his original monograph \cite{Lyap1892}, Lyapunov studied ordinary differential equations and provided two
results in particular that would have wide-ranging impact:
\begin{theorem}{\cite[Section 16, Theorem I]{Lyap1892}}
	\label{thm:OriginalLyapunov}
	If the differential equations of the disturbed motion are such that it is possible to find a
	definite function $V$, of which the derivative $V'$ is a function of fixed sign which is opposite
	to that of $V$, or reduces identically to zero, the undisturbed motion is stable.
\end{theorem}

\begin{theorem}{\cite[Section 16, Remark II]{Lyap1892}}
	\label{thm:LyapunovRemark}
	If the function $V$, while satisfying the conditions of the theorem, admits an infinitely small upper limit,
	and if its derivative represents a definite function, we can show that every disturbed motion, sufficiently
	near the undisturbed motion, approaches it asymptotically.
\end{theorem}

For simplicity, ``undisturbed motion'' can be taken to be an equilibrium point of a differential equation while
``disturbed motion'' refers to solutions of a differential equation originating from a point other than
the equilibrium point; i.e., motions that are {\rev initially disturbed or} perturbed away from the equilibrium.
{\rev However, any non-trivial solution of an ordinary differential equation can be considered where a
desired ``reference'' solution is the undisturbed motion and the disturbed motion refers to solutions
perturbed away from the reference solution.  Similarly, one may consider general attractors whereby a
disturbed motion is one that originates outside the attractor.}

The strength of Lyapunov's second method as encapsulated in Theorems \ref{thm:OriginalLyapunov} and
\ref{thm:LyapunovRemark} is that it is possible
to ascertain stability without solving the underlying differential equation.  However, the difficulty of
Theorems \ref{thm:OriginalLyapunov} and \ref{thm:LyapunovRemark} lies in finding an appropriate
function $V$.
Therefore, the converse or existence question arises; i.e., if an
``undisturbed motion'' or equilibrium point is stable or asymptotically stable, does an appropriate
function $V$ exist?  A related question is: how can such a function be constructed?  The first question
is the subject of this survey paper, while the second is very much the subject of ongoing research.

In addition to their intrinsic mathematical interest, converse Lyapunov theorems are important in
that they indicate which stability properties can always be established by an appropriate Lyapunov function.
In fact, the study of the converse question was crucial in discovering that Theorem \ref{thm:LyapunovRemark}
implies more than asymptotic stability and, in fact, implies {\it uniform} {\rev (with respect to initial time)} asymptotic stability.
Converse Lyapunov theorems are also a useful tool when considering perturbed systems, where the perturbations
may be additive to the system equations \cite{Barb50-DAN}, \cite{Malk54-PMM}, time-delays in the system equations \cite{Kras63}, or
as the result of a linearization \cite{Lyap1892}.  It was, in fact, this latter concern that motivated Lyapunov to
develop his second method.  Finally, in the development of numerical, constructive techniques for Lyapunov functions,
converse Lyapunov theorems provide a gold-standard by which these techniques can be measured\footnote{This
is similar to the role played by Shannon's channel coding theorem in information theory \cite{ShWe49}
(see also \cite[Section 7.7]{CoTh06}).
Shannon's theorem provides
a fundamental limit for communication over a noisy channel by showing that a capacity-achieving channel coding scheme must exist,
but does not constructively provide a coding scheme that
achieves that limit.  Nonetheless, Shannon's limit has been an invaluable idealized goal for information and coding theorists for over
60 years.}.  For a thorough survey of computational methods for Lyapunov functions, see \cite{GiHa15-DCDSB},
which is in the same special issue as this paper.

This survey is organized as follows.  In Section \ref{sec:LyapunovSecondMethod} we
provide the basic theory of Lyapunov's second method.  In Section \ref{sec:constructive} we describe how
Lyapunov functions can be constructed for linear systems via an algebraic approach and in {\rev Section \ref{sec:Lure}
we discuss extensions to linear systems in feedback with static nonlinearities; i.e., so-called Lur'e systems.}
In Section \ref{sec:Zubov} we present
a constructive technique for autonomous systems based on solution of a partial differential equation.  In the results of
both Section \ref{sec:constructive} and Section \ref{sec:Zubov},
the assumption of asymptotic stability guarantees that the described techniques yield a Lyapunov function.
In Section \ref{sec:historical} we trace the historical development of converse Lyapunov theorems and briefly
describe some of the approaches used.  In Section \ref{sec:KellettTeel} we present some specific converse
Lyapunov theorems for {\rev so-called $\KL$-stability of} differential and difference inclusions; {\rev a concept
equivalent to uniform global asymptotic stability}.  In Section \ref{sec:Instability} we present results for unstable equilibrium points.  Some concluding remarks are contained in Section \ref{sec:Conclusions}.

% -----------------------------------------------------------------------------------------------------------------------------------------
%		LYAPUNOV'S SECOND METHOD
% -----------------------------------------------------------------------------------------------------------------------------------------

\section{Lyapunov's second method}
\label{sec:LyapunovSecondMethod}

Initially, we will consider dynamical systems described by ordinary differential equations
\begin{equation}
	\label{eq:sys}
	\dot{x} = f(x,t)
\end{equation}
where, for simplicity, we assume that $f: \R^n \times \R_{\geq 0} \rightarrow \R^n$ is locally Lipschitz in
$x$ and continuous in $t$ so that (local)
existence and uniqueness of solutions is guaranteed.  {\rev Forward completeness will generally follow from
assumed stability properties or the presence of a Lyapunov function.}  We denote a solution to \eqref{eq:sys}
from the initial state $x \in \R^n$ and initial time $t_0 \in \R_{\geq 0}$ at the time $t \geq t_0$ by
$\phi(t,t_0,x)$.  In other words, $\phi: \R_{\geq 0} \times \R_{\geq 0} \times \R^n \rightarrow \R^n$ satisfies
\begin{equation}
	\frac{d}{dt} \phi(t,t_0,x) = f(\phi(t,t_0,x),t),
\end{equation}
and we write $\phi \in \S_{t_0}(x)$ where $\S_{t_0}(x)$ denotes the set of solutions from initial
time $t_0 \in \R_{\geq 0}$ and initial state $x \in \R^n$.  As we initially assume uniqueness of solutions
for \eqref{eq:sys}, $\S_{t_0}(x)$ contains a single function.  In the sequel when we consider
difference or differential inclusions, or \eqref{eq:sys} where the righthand side is only continuous,
the set $\S_{t_0}(x)$ will, in general, be larger than merely
a singleton.
We further assume that $f(0,t) = 0$ for all $t \geq t_0$ so that the origin is an equilibrium point.

Lyapunov precisely defined the notion of stability.  However, many types of stability are possible
in general and the most useful notions, presented below, are largely due to Chetaev \cite{Chet61},
Malkin \cite{Malk54-PMM}, Massera \cite{Mass56-AM}, and Barbashin and Krasovskii \cite{BaKr52-DAN}.

Throughout this survey we will make use of the comparison functions introduced
by Massera \cite{Mass56-AM} and Hahn \cite{Hahn67}.  The use of such functions simplifies
many statements and proofs in the area of systems theory.  {\rev To denote the class of positive definite
functions with domain $\R_{\geq 0}$ we use the notation $\mathcal{P}$ ($\rho \in \mathcal{P}$); i.e.,
functions $\rho: \R_{\geq 0} \rightarrow \R_{\geq 0}$ that are continuous, zero at zero, and strictly
positive elsewhere.}  A function $\alpha: \R_{\geq 0}
\rightarrow \R_{\geq 0}$ is said to be of class-$\K$ ($\alpha \in \K$) if it is continuous, zero at
zero, and strictly increasing.  It is said to be of class-$\Kinf$ ($\alpha \in \Kinf$) if, in addition, it approaches infinity
as its argument approaches infinity.  A function $\sigma : \R_{\geq 0} \rightarrow \R_{\geq 0}$ is
said to be of class-$\mathcal{L}$ ($\sigma \in \mathcal{L}$) if it is continuous, strictly decreasing, and
approaches zero as its argument
approaches infinity.  A function $\beta : \R_{\geq 0} \times \R_{\geq 0} \rightarrow \R_{\geq 0}$ is
said to be of class-$\KL$ ($\beta \in \KL$) if it is of class-$\K$ in its first argument and of class-$\mathcal{L}$
in its second argument.
For a more extensive introduction to such comparison functions, see \cite{Kell14-MCSS}.

In what follows, for a set $\D \subset \R^n$ we denote its boundary by $\partial \D$ and its closure by
$\overline{\D}$.  We denote the open ball of radius $\varepsilon \in \R_{> 0}$, {\rev centered at the origin,} by
\[ \B_\varepsilon := \left\{ x \in \R^n : |x| < \varepsilon \right\} \]
and we write $\B = \B_1$.
\begin{definition}
	The origin is said to be {\em stable for \eqref{eq:sys}} if {\rev there exists a neighborhood of the origin
	$\mathcal{N} \subset \R^n$ so that} for each $t_0 \in \R_{\geq 0}$ there exists
	$\alpha_{t_0} \in \K$ so that, for all $x \in \mathcal{N}$ and
	all $t \geq t_0$,
	\begin{equation}
		|\phi(t,t_0,x)| \leq \alpha_{t_0}(|x|).
	\end{equation}
 The origin is said to be {\em uniformly stable for \eqref{eq:sys}} if the
	function $\alpha_{t_0} = \alpha \in \K$ can be chosen independent of the initial time $t_0 \in \R_{\geq 0}$.
\end{definition}

\begin{definition}
	The origin is said to be {\em asymptotically stable for \eqref{eq:sys}} if {\rev there exists a neighborhood of the origin
	$\mathcal{N} \subset \R^n$ so that} for each $x \in \mathcal{N}$ and $t_0 \in \R$
	there exists {\rev $\sigma_{x,t_0} \in \L$
	so that, for all $t \geq t_{0}$,
	\begin{equation}
		\label{eq:as_stable}
		|\phi(t,t_0,x)| \leq \sigma_{x,t_0}(t-t_0).
	\end{equation}}The origin is said to be {\em equiasymptotically stable for \eqref{eq:sys}} if, {\rev for every $t_0 \in \R$,
	there exists a function $\beta_{t_0} \in \KL$ so that, for all $t \geq t_0$,
	\begin{equation}
		|\phi(t,t_0,x)| \leq \beta_{t_0}(|x|,t-t_0), \quad \forall x \in \mathcal{N}.
	\end{equation}}The origin is said to be {\em uniformly asymptotically stable for \eqref{eq:sys}} if the function $\beta_{t_0} = \beta \in \KL$
	can be chosen independent of the initial time
	$t_0 \in \R_{\geq 0}$.
\end{definition}

{\rev To clarify, the difference between asymptotic stability and equiasymptotic stability lies in the fact that the latter
is uniform with respect to the size of the initial state, while uniform asymptotic stability requires that
the stability property is uniform with respect to both the size of the initial state and the initial time.
One can also define the property of {\it uniformly attractive} where asymptotic stability is uniform
with respect to the initial time but not the initial state; i.e., $\sigma_{x,t_0} = \sigma_x \in \L$ in \eqref{eq:as_stable}.
However, this property appears to have found limited use and we do not consider it further.
}

{\rev Note that in the above definitions, the existence condition for the neighborhood $\mathcal{N} \subset \R^n$ leads to these
stability properties sometimes being refered to as {\em local}.  By contrast, for an open set $\G \subset \R^n$
fixed {\em a priori} and} containing the origin, the above stability properties are said to hold ``in the large'' if they
hold for all $x \in \G$.  {\rev For example, the origin is said to be asymptotically stable in the large on $\G \subset \R^n$
if for every $x \in \G$ and $t_0 \in \R$ there exists $\sigma_{t_0} \in \L$ so that \eqref{eq:as_stable} holds.}
In the event that $\G = \R^n$, the above stability properties are said to be ``global''.

{\rev Massera \cite[Theorem 7]{Mass56-AM} provided several relationships amongst the above stability properties.
Of particular interest, uniform stability implies stability and uniform asymptotic stability implies equiasymptotic stability
which, in turn, implies asymptotic stability.  Furthermore, \cite[Theorem 7]{Mass56-AM} demonstrates that if the righthand
side of \eqref{eq:sys} is periodic in, or independent of, the time $t$, then the converses hold; i.e., stability implies uniform
stability and asymptotic stability implies uniform asymptotic stability which, in turn, implies equiasymptotic stability.}

The classical definitions for the various stability concepts are given in $\varepsilon$-$\delta$ terms for stability
and as a combination of stability and limiting behavior as time approaches infinity for the asymptotic stability
concepts.  That these definitions are equivalent to the comparison function formulations presented above
was shown by Hahn \cite{Hahn67}.

The following theorem summarizes Lyapunov's second method as it relates to \eqref{eq:sys} and the stability definitions
described above; see \cite{Hahn67}, \cite{Kras63}, or \cite{RHL77}.
\begin{theorem}
	\label{thm:ModernLyapunov}
	Let $V: \R^n \times \R \rightarrow \R_{\geq 0}$ be a continuously differentiable function and
	consider the following conditions on $V$:
	\begin{enumerate}[(i)]
		\item \label{it:pos_def}
			Suppose there exists $\alpha_1 \in \K$ so that, for all $x \in \R^n$ and $t \geq t_0$
			\begin{equation}
				\label{eq:V_lower}
				\alpha_1(|x|) \leq V(x,t).
			\end{equation}
{\rev		\item \label{it:strong_pos_def}
			Suppose there exists a continuous function $\kappa: \R_{\geq 0} \times \R_{\geq 0} \rightarrow \R_{\geq 0}$
			such that $\kappa(\cdot,t) \in \K$ for fixed $t \in \R$ and $\kappa(s,\cdot)$ is continuous, positive,
			monotone increasing, and unbounded for each fixed $s \in \R_{\geq 0}$.  Furthermore, suppose
			that, for all $x \in \R^n$ and $t \geq t_0$,
			\begin{equation}
				\label{eq:V_strong_pos_def}
				\kappa(|x|,t) \leq V(x,t).
			\end{equation} }
		\item \label{it:decrease}
		 	Suppose there exists $\rho : \R_{\geq 0} \rightarrow \R_{\geq 0}$ continuous and, for all
			$x \in \R^n$ and $t \geq t_0$
			\begin{equation}
				\label{eq:decrease}
				\left. \tfrac{d}{dt} V(\phi(t,t_0,x),t) \right|_{t = t_0} = \tfrac{\partial}{\partial t} V
				+ \left\langle \tfrac{\partial}{\partial x}V, f(x,t) \right\rangle
				\leq - \rho(|x|) .
			\end{equation}
		\item \label{it:decrescent}
			Suppose there exists $\alpha_2 \in \K$ so that, for all $x \in \R^n$ and $t \geq t_0$
			\begin{equation}
				\label{eq:V_upper}
				V(x,t) \leq \alpha_2(|x|) .
			\end{equation}
	\end{enumerate}
	The following statements hold with respect to \eqref{eq:sys}:
	\begin{enumerate}[(a)]
		\item \label{it:stable} If items (\ref{it:pos_def}) and (\ref{it:decrease}) hold, then the origin is stable;
		\item If items (\ref{it:pos_def}), (\ref{it:decrease}), and (\ref{it:decrescent}) hold, then the
			origin is uniformly stable;
		\item If items (\ref{it:strong_pos_def}) and (\ref{it:decrease}) hold, then the origin
			is equiasymptotically stable;
		\item If items (\ref{it:pos_def}), (\ref{it:decrease}), and (\ref{it:decrescent}) hold with
			$\rho \in \mathcal{P}$, then the origin is uniformly asymptotically stable;
		\item \label{it:UGAS} If items (\ref{it:pos_def}), (\ref{it:decrease}), and (\ref{it:decrescent}) hold with
			$\rho \in \mathcal{P}$ and
			$\alpha_1, \alpha_2 \in \Kinf$, then the origin
			is uniformly globally asymptotically stable.
	\end{enumerate}	
\end{theorem}

We will refer to \eqref{eq:decrease} as the ``derivative of $V$ along solutions of \eqref{eq:sys}'' or
simply as the ``total derivative of $V$''.
We note that continuous differentiability of the Lyapunov function $V$ is not critical to the development of
the theory as the key idea is that the Lyapunov function should decrease along solutions of \eqref{eq:sys}.
This property can be stated without any requirements on the regularity of $V$, but then may require explicit
knowledge of the solutions.  While continuous differentiability leads to the simple criterion of \eqref{eq:decrease},
decrease conditions involving nonsmooth derivatives (such as Dini derivatives or subgradients) for functions with
weaker regularity properties can be used (see, e.g., \cite{CLSW98}, \cite{MHL08}, \cite{Roxi65a-JDE}, and the references therein).

\begin{remark}
	\label{rem:InTheLarge}
	A sufficient condition for (asymptotic) stability in the large on an open set $\G \subset \R^n$ containing the origin
	in the above theorem is the requirement that $\lim_{x \rightarrow \partial \G} V(x) = \infty$, where, in directions in
	which $\G$ is unbounded this is interpreted as $\lim_{|x| \rightarrow \infty} V(x) = \infty$.  Define
	\[ \omega(x) := \max \left\{ |x|, \frac{1}{|x|_{\partial \G}} - \frac{1}{|0|_{\partial \G}} \right\} \]
	where $|x|_{\partial \G} := \min_{y \in \partial \G} |x-y|$ denotes the closest distance to the boundary of $\G$.
	Then, to guarantee in-the-large stability properties, we
	can replace the lower bound \eqref{eq:V_lower} by
	\[ \alpha_1(\omega(x)) \leq V(x,t), \quad \forall x \in \G, \ t \in \R_{\geq 0}  \]
	and with the requirement that $\alpha_1 \in \Kinf$.
	If $\G = \R^n$ then this implies $\omega(x) = |x|$ and stability in the large coincides with global stability.
	This property of $V$ is refered to as {\em radially unbounded} (on $\G$) by Hahn \cite{Hahn67} and as
	{\em infinitely large} by Barbashin and Krasovskii \cite{BaKr52-DAN}.
\end{remark}

The property of $V$ described in item \eqref{it:decrescent} was termed {\em decrescent} by Hahn \cite{Hahn67}
and as an infinitely small upper bound by Lyapunov \cite{Lyap1892} as in Theorem \ref{thm:LyapunovRemark},
where the fact that the upper bound is infinitely small clearly only holds near the origin.  Note that if
$V$ is continuous, independent of $t$, and satisfies $V(0) = 0$, then the decrescent bounds hold trivially.
Hahn refered to the property described in Theorem \ref{thm:ModernLyapunov}.\ref{it:strong_pos_def} as strongly positive
definite \cite[Definition 41.5]{Hahn67}.

The converse question, then, is which of the statements in Theorem \ref{thm:ModernLyapunov} can be
reversed?  For systems described by ordinary differential equations and for stability properties related
to the origin, this question was largely answered by the end of the 1950's (see Section \ref{sec:historical}
below).  For more general systems and more general stability properties, research is still ongoing
(see Section \ref{sec:KellettTeel} below).  Prior to addressing the general existence result, we discuss three
constructive techniques.

% -----------------------------------------------------------------------------------------------------------------------------------------
%		CONSTRUCTIVE METHODS / LINEAR SYSTEMS
% -----------------------------------------------------------------------------------------------------------------------------------------

\section{Linear systems}
\label{sec:constructive}
For general systems such as \eqref{eq:sys}, finding an explicit closed-form Lyapunov function
is known to be a difficult task.  However, for linear systems, finding such a Lyapunov function
is essentially an algebraic problem.  {\rev In what follows, we denote a positive definite symmetric
matrix $P \in \R^{n \times n}$ by $P > 0$.}

The first converse theorem was demonstrated in Lyapunov's original monograph
for the case of linear systems
\begin{equation}
	\label{eq:CTlinsys}
	\dot{x} = Ax, \quad x \in \R^n .
\end{equation}
Stated in modern terms, the following is \cite[Section 20, Theorem II]{Lyap1892}:
\begin{theorem}
	Given any $Q>0$ there exists $P>0$ satisfying
	\begin{equation}
		\label{eq:LyapunovEqn}
		A^TP + PA = -Q
	\end{equation}
	if and only if the origin is asymptotically stable for \eqref{eq:CTlinsys}.
\end{theorem}

Equation \eqref{eq:LyapunovEqn} is refered to as the {\em Lyapunov equation} and straightforward calculations
show that the quadratic $V(x) = x^TPx$ is a Lyapunov function for \eqref{eq:CTlinsys}.

A similar result holds for linear discrete time systems described by
\begin{equation}
	\label{eq:DTlinsys}
	x^+ = Ax, \quad x \in \R^n .
\end{equation}
The following result is due to Stein \cite[Theorem 1]{Stei52-JRNBS} and was first mentioned in
a systems theoretic context in \cite{Hahn63} and \cite{KaBe60b-JBE}.
\begin{theorem}
	Given any $Q>0$ there exists $P>0$ satisfying
	\begin{equation}
		\label{eq:DT_LyapunovEqn}
		A^TPA - P = -Q
	\end{equation}
	if and only if the origin is asymptotically stable for \eqref{eq:DTlinsys}.
\end{theorem}
As in the continuous time case, \eqref{eq:DT_LyapunovEqn} is called the {\em discrete time Lyapunov
equation} (or the Stein equation after \cite{Stei52-JRNBS}) and the quadratic $V(x) = x^TPx$ is a Lyapunov function for \eqref{eq:DTlinsys}.

It has proved difficult to generate converse theorems by directly constructing Lyapunov functions
for systems more general than linear time-invariant systems.  This can be seen in the conditions available
for linear time-varying systems where \eqref{eq:LyapunovEqn} or \eqref{eq:DT_LyapunovEqn}
are replaced by matrix differential or difference equations \cite[Theorem 5]{AnMo69-SJC},
\cite[Theorem 4.3]{AnMo81-SICON}, for which closed form solutions generally do not exist.
General nonlinear systems, naturally, present an even greater challenge.

% -----------------------------------------------------------------------------------------------------------------------------------------
%		CONSTRUCTIVE METHODS / LUR'E SYSTEMS
% -----------------------------------------------------------------------------------------------------------------------------------------

{\rev
\section{Lur'e systems}
\label{sec:Lure}
Due to its importance in engineering applications, significant efforts were made to extend the
converse theorems for linear systems described above to the case of so-called
{\em absolute stability}; i.e., for asymptotic stability of the origin for linear systems
with a static, memoryless, sector-bounded nonlinearity in a feedback loop.
Such systems have come to be called {\em Lur'e}\footnote{The Cyrillic-to-Latin transliteration
of {\cyr Lur{\cprime}e} has led to four distinct spellings in the western literature: Lure, Lur'e, Lurie, and Lur\'{e}.}
{\em systems} and can be written as
\begin{equation}
	\label{eq:Lure}
	\begin{array}{rcl}
		\dot{x} & = & Ax - B\psi(y) \\
		y & = & Cx
	\end{array}
\end{equation}
where $x \in \R^n$, $y \in \R^p$, and $\psi: \R^p \rightarrow \R^p$ satisfies the sector bound
\begin{equation}
	\label{eq:sector_cdn}
	\langle \psi(y), y - K\psi(y) \rangle \geq 0
\end{equation}
for some positive definite symmetric matrix $K \in \R^{p \times p}$.  Note that, with the above
definitions, the system is assumed to have the same number of inputs as outputs; i.e.,
$B \in \R^{n \times p}$ and $C \in \R^{p \times n}$.

The study of such systems was instigated by Lur'e and Postnikov \cite{LuPo44-PMM}
with the motivation of studying the stability of controlled systems subject to common
nonlinear actuators.  In particular, the matrix $A$ is assumed to be Hurwitz (i.e., to have
all its eigenvalues in the open left half plane), where this property may have been imposed
by a (linear) feedback control, and the unknown but sector bounded nonlinearity may correspond
to the feedback being implemented by an actuator with unmodeled characteristics such as friction
or deadzones.

Despite the importance of such systems, and the significant effort expended in trying to derive
a converse theorem for absolute stability, a constructive converse has yet to be found.
However, for single-input single-output systems, an existence result based on a frequency domain condition
was provided in three fundamental papers by Popov \cite{Popo63-AiT},
Yakubovich \cite{Yaku62-DAN}, and Kalman \cite{Kalm63-PNAS}.  The modern
form of this result is called the Popov-Yakubovich-Kalman Lemma or the Positive Real Lemma (e.g., \cite[Lemma D.6]{KKK95}
or the thorough discussion in \cite[Appendix H]{Silj69}).  It is beyond the scope
of this survey to deal with the Popov-Yakubovich-Kalman Lemma and its significant ramifications\footnote{The
interested reader is directed to \cite[Appendix H]{Silj69}, \cite[Section 10.1]{Khal96}, \cite{KoAr01-Aut}, and \cite{SWMWK07-SIAM-Review}.}.  The extension from single-input single-output systems to multi-input multi-output systems was provided
independently by Popov \cite[Lemma 9.3.1]{Popo66} (also \cite{Popo64}) and Anderson \cite{Ande66-JFI}.

\begin{theorem}
	\label{thm:Lure}
	Suppose that, for \eqref{eq:Lure},  $A$ is Hurwitz, $(A,B)$ is controllable, $(C,A)$ is
	observable.  Let $G(s) := C(sI-A)^{-1}B$.  Suppose the positive definite symmetric matrix
	$K$ is such that $K\psi(y)$ is the gradient of a positive semidefinite scalar function; i.e.,
	\begin{equation}
		\label{eq:Ksemidefinite}
		\int_{\Gamma(0,y)} \psi^T(s)K ds \geq 0, \quad \forall y \in \R^p,
	\end{equation}
	where $\Gamma$ is any smooth curve in $\R^p$ connecting $0$ and $y = Cx$.
	For any $\eta \in \R_{>0}$ such that $1/\eta$ is not an eigenvalue of $A$, if $Z(s) := I + (1+\eta s)KG(s)$
	is strictly positive real, then there exists a positive definite symmetric matrix $P$ so that
	\begin{equation}
		\label{eq:Lure-Postnikov}
		V(x) = x^T P x + \int_{\Gamma(0,y)} \psi^T(s)K ds
	\end{equation}
	is a Lyapunov function.
\end{theorem}

Note that since $K\psi(y)$ is the gradient of a scalar function, the integral \eqref{eq:Ksemidefinite} is path-independent
\cite[Theorem 10-37]{Apos57}.

As in \cite{Ande66-JFI}, an example $K$ and $\psi$ satisfying the conditions of Theorem~\ref{thm:Lure} is
when $K$ is diagonal and the sector-bounded nonlinearities are decoupled so that each of $p$ nonlinearities only depends on one
element of $y$.  Other examples are provided in \cite[Section 10.1]{Khal96}.

We observe that Theorem~\ref{thm:Lure} is not a converse theorem in the usual sense in that it does not
start from a stability property and then provide a Lyapunov function.  Rather, it starts from the
requirement that the transfer function matrix $Z(s)$ be strictly positive real (see
\cite[Definition 10.3]{Khal96} for a definition of strict positive realness), which by the Popov-Yakubovich-Kalman
Lemma implies solvability of the matrix equations
\begin{equation}
	\label{eq:PR_cdn}
	\begin{array}{rl}
	P\hat{A} + \hat{A}^TP & = -L^TL - \varepsilon P \\
	P\hat{B} & = \hat{C}^T - L^TW \\
	W^TW & = \hat{D} + \hat{D}^T
	\end{array}
\end{equation}
for matrix $L$, positive definite symmetric matrix $P$, and constant $\varepsilon \in \R_{>0}$, where
\[ \hat{A} := A, \quad \hat{B} := B, \quad \hat{C} := KC + \eta KCA, \quad \hat{D} := I + \eta KCB. \]
The matrix $P$ from \eqref{eq:PR_cdn} is then the $P$ of the Lur'e-Postnikov Lyapunov function
\eqref{eq:Lure-Postnikov}.  The lack of a standard converse
result highlights the difficulty inherent in finding Lyapunov functions for general nonlinear
systems.

Similar results on absolute stability and Lyapunov functions are available for discrete time systems
with the original study of such systems performed in a sampled-data context by Tsypkin
\cite{Tsyp62-Doklady}, \cite{Tsyp63-AiT}.  The discrete time version of the matrix conditions \eqref{eq:PR_cdn}
was provided independently by Popov \cite[Theorem 10.1.1]{Popo66} and Hitz and Anderson \cite{HiAn69-IEE} and is similar to the difference between the Lyapunov equation \eqref{eq:LyapunovEqn} and the Stein equation \eqref{eq:DT_LyapunovEqn}:
\begin{equation}
%	\label{eq:PR_cdn}
	\begin{array}{rl}
	\hat{A}^T P\hat{A} - P & = -L^TL  \\
	\hat{A}^T P\hat{B} & = \hat{C}^T - L^TW \\
	W^TW & = \hat{D} + \hat{D}^T - \hat{B}^TP \hat{B} .
	\end{array} \nonumber
\end{equation}

The study of controlled systems, and systems subject to external disturbances, that is indicated by the
structure of \eqref{eq:Lure}, was generalized to the notion of dissipativity \cite{Will-ARMA72a, Will-ARMA72b}
and later to input-to-state stability \cite{Sont89-TAC, SoWa95-SCL}, {\rev input-output-to-state stability
\cite{KSW01-SICON} and measurement-to-error stability \cite{ISW02-CDC}}.  While Lyapunov methods are
of significant importance in these topics, and existence results are available, the study of such systems
is beyond the scope of this survey.

}

% -----------------------------------------------------------------------------------------------------------------------------------------
%		ZUBOV'S METHOD
% -----------------------------------------------------------------------------------------------------------------------------------------

\section{Zubov's method}
\label{sec:Zubov}

Zubov \cite{Zubo64} presented a method for estimating the domain of attraction of an autonomous
ordinary differential equation
\begin{equation}
	\label{eq:autonomous}
	\dot{x} = f(x)
\end{equation}
where $f:\R^n \rightarrow \R^n$ is locally Lipschitz and $f(0) = 0$.  In particular, Zubov's method
constructs a Lyapunov function that guarantees asymptotic stability in the large on the domain of attraction.

If the origin is uniformly asymptotically stable for \eqref{eq:autonomous}, then the domain of attraction $\D \subset \R^n$
for  the origin is
\begin{equation}
	\D := \left\{ x \in \R^n : \lim_{t \rightarrow \infty} \phi(t,x) = 0 \right\}.
\end{equation}
We observe that $\D$ is an open set.  {\rev For $\G \subset \R^n$, a function $V: \G \rightarrow \R$ is positive definite if
$V(0) = 0$ and $V(x) > 0$ for all $x \in \G \backslash \{0\}$.}

{\rev The following theorem combines \cite[Theorem 34.1]{Hahn67} and \cite[Theorem 51.1]{Hahn67}.}
\begin{theorem}
	\label{thm:Zubov}
	The origin is asymptotically stable for \eqref{eq:autonomous} on a domain of attraction $\D \subset \R^n$
	if and only if there exist functions $V,h: \R^n \rightarrow \R_{\geq 0}$ satisfying
	\begin{enumerate}[(i)]
	\item $V$ is continuous and positive definite in $\D$, $0 \leq V(x) < 1$, and
		$\lim_{|x| \rightarrow \partial \D} V(x)\break = 1$;
	\item $h$ is continuous {\rev and positive definite}; and
	\item the following partial differential equation is satisfied
	\begin{equation}
		\label{eq:Zubov_PDE}
		\left\langle \tfrac{\partial}{\partial x}V(x), f(x) \right\rangle = - h(x)(1-V(x)) \sqrt{1+ |f(x)|^2}.
	\end{equation}
	\end{enumerate}
\end{theorem}

In other words, given a uniformly asymptotically stable equilibrium point, it is always possible to find a Lyapunov
function, defined on the domain of attraction, that satisfies the partial differential equation \eqref{eq:Zubov_PDE}.

To construct an appropriate Lyapunov function, we start from the characterization of uniform
asymptotic stability given by $\alpha \in \Kinf$, $\sigma \in \mathcal{L}$, and
\begin{equation}
	|\phi(t,x)| \leq \alpha(|x|) \sigma(t), \quad \forall x \in \D, \ t \in \R_{\geq 0}.
\end{equation}
{\rev Define
\[ \sigma^\dag(s) := \left\{ \begin{array}{ccl}
	\sigma^{-1}(s) & , & s \in (0,\sigma(0)] \\
	0 & , & s \geq \sigma(0)
\end{array} \right. \]
and
\begin{equation}
	\label{eq:scaling}
	\varphi(s) := \left\{ \begin{array}{ccl}
		s \exp(-\sigma^{\dag}(s)) & , & s > 0 \\
		 0 & , & s = 0 .
		\end{array} \right.
\end{equation}}That $\varphi \in \Kinf$ follows from basic compositional properties
of comparison functions; see \cite{Kell14-MCSS}.  The functions
\begin{align}
	\label{eq:Zubov}
	V(x) & := 1 - \exp\left(- \int_0^\infty \varphi(|\phi(\tau,x)|)d\tau\right), \quad {\rm and} \\
	h(x) & := \frac{\varphi(|x|)}{\sqrt{1 + |f(x)|^2}}
\end{align}
can then be shown to satisfy the necessary properties of Theorem~\ref{thm:Zubov}.  See
\cite[Theorem 51.1]{Hahn67} for a complete proof.

We observe that, by definition, the origin is uniformly asymptotically stable in the large on its domain
of attraction.  In Remark \ref{rem:InTheLarge} we observed that a Lyapunov function that is radially
unbounded on $\D$ can be used to conclude in-the-large stability properties.  However, as we see in
Zubov's method, the derived Lyapunov function approaches the value $1$ on the boundary of the
domain of attraction.  The critical observation is that this leads to $V$ being {\em proper} on the domain
of attraction, where the term proper refers to preimages of compact sets being compact.  In other words,
for any compact set $[0,c] \subset [0,1)$, the set defined by $V^{-1}([0,c])$ must also be
compact in $\R^n$.  If $V$ is proper and the time derivative of $V$ along solutions is negative definite, then
trajectories necessarily move from larger (compact) level sets to smaller (compact) level sets.  The radially unbounded on
$\D$ condition of Remark \ref{rem:InTheLarge} is a sufficient condition for $V$ to be proper.\footnote{ {\rev
In order to differentiate between a definition of {\em proper} that requires that the inverse of all compact sets
to be compact (i.e., $[0,c] \subset \R_{\geq 0}$) and one that requires that the inverse of compact sets
on the range of the function $V$ to be compact (i.e., $[0,c] \subset [0,1)$ as above), the terminology
{\em proper on its range} or {\em semiproper} is sometimes used.}}

Furthermore, we note that $\mu(s) := - \log(1-s)$ maps $[0,1) \mapsto [0,\infty)$, is continuously
differentiable, and strictly increasing.  Consequently, the function $W(x) := \mu(V(x))$ will be a Lyapunov
function for asymptotic stability of the origin in the large on $\D$ as defined by Barbashin and Krasovskii where
$\lim_{x \rightarrow \partial \D} W(x) = \infty$.  Straightforward manipulations of \eqref{eq:Zubov} yield that
\begin{equation}
	\label{eq:semi_infinite_integral_Zubov}
	W(x) = \int_0^\infty \varphi(|\phi(\tau,x)|) d\tau.
\end{equation}
We will see this Lyapunov function candidate again in \eqref{eq:Massera} below.

Zubov \cite[Theorems 19 and 78]{Zubo64} extended the above result to dynamical systems on metric spaces
including time-varying systems and systems described by PDEs that admit classical solutions, as well
as accounting for asymptotic stability of closed invariant sets as opposed to merely the origin.

While Zubov's method is not constructive in the same sense that solving the Lyapunov equation
\eqref{eq:LyapunovEqn} is constructive, the freedom to choose the function $h$ in \eqref{eq:Zubov_PDE}
has enabled useful constructions in many cases.
See \cite{GCW01-SICON}, \cite{GrSe11-SICON} and the references therein for recent applications
and extensions of Zubov's method.

% -----------------------------------------------------------------------------------------------------------------------------------------
%		LYAPUNOV FUNCTION CONSTRUCTION IN CONVERSE THEOREMS
% -----------------------------------------------------------------------------------------------------------------------------------------

\section{Historical developments}
\label{sec:historical}

Though the converse question was answered by Lyapunov in the linear case, the converse of
Lyapunov's second method in the more general
case represented by \eqref{eq:sys} remained open throughout the early 1900's.

\subsection{Early results - pre-1950}
Persidskii \cite{Pers37-Doklady} provided the first general converse theorem when he demonstrated
that, under the assumption that the origin is a stable equilibrium point, the function
\begin{equation}
	V(x,t) = \min_{t_0 \leq \tau \leq t} |\phi(\tau,t,x)|
\end{equation}
is in fact a Lyapunov function.

In \cite{Mass49-AM}, Massera precisely defined stability, asymptotic stability, and
equiasymptotic stability and compared them via examples as well as by sufficient Lyapunov
function properties that guaranteed them.  In the case when \eqref{eq:sys} is periodic or autonomous,
and the origin is asymptotically stable, Massera showed that the semi-infinite integral
\begin{equation}
	\label{eq:Massera}
	V(x,t) = \int_t^\infty \alpha(|\phi(\tau,t,x)|)  d\tau
\end{equation}
where $\alpha: \R_{\geq 0} \rightarrow \R_{\geq 0}$ is an appropriately chosen
continuous function, is in fact a continuously differentiable Lyapunov function.
Furthermore, Massera demonstrated that if \eqref{eq:sys} is periodic in $t$ or independent of
$t$ then $V$ has this same property.

Massera's manuscript \cite{Mass49-AM} would have a significant impact on the study of the converse
question.  Not only did Massera provide the first converse theorem for asymptotic stability, but
\cite{Mass49-AM} left open the converse question for systems that were neither periodic nor autonomous
(a problem that required the notion of uniform stability as described below).  The question of the
existence of a {\it smooth} Lyapunov function remained, as did whether or not the assumption in
\cite{Mass49-AM} of continuous differentiability of the righthand side of \eqref{eq:sys} was necessary.
Finally, the proof technique used by Massera became the standard approach in much subsequent work.
In particular, most subsequent authors have proposed Lyapunov function candidates similar to the semi-infinite
integral of \eqref{eq:Massera} and, frequently, the choice of the scaling $\alpha$ is done either directly from,
or similar to, that from what is now
called ``Massera's Lemma'' \cite[p.\ 716]{Mass49-AM}.

Contemporaneously with \cite{Mass49-AM}, and using a different proof technique, Barbashin \cite{Barb50-DAN}
demonstrated that, for an autonomous system \eqref{eq:sys}, there exists a Lyapunov
function with the same regularity as that of the vector field $f$.

\subsection{Fundamental theory - 1950's}
Malkin \cite{Malk54-PMM} recognized that the important, and more general, property that
allowed Massera to derive converse theorems for periodic and autonomous systems
with an asymptotically stable equilibrium point is that of uniformity with respect to time.
Furthermore, Malkin demonstrated that, as it was originally written in \cite[Section 16, Remark 2]{Lyap1892}
(Theorem \ref{thm:LyapunovRemark} in Section \ref{sec:Intro}),
Lyapunov's second method in fact guarantees {\it uniform} asymptotic stability.  In particular, the uniformity
follows from the assumed decrescent property of the Lyapunov function

Around the same time, Barbashin and Krasovskii \cite{BaKr52-DAN} demonstrated that a sufficient condition
for (asymptotic) stability in the large on a set $\G \subseteq \R^n$ is that the Lyapunov function be radially unbounded
on $\G$.  Subsequently, following Malkin and Massera's proof techniques,
Barbashin and Krasovskii \cite{BaKr54-PMM} demonstrated that a radially unbounded (on $\G$) and decrescent Lyapunov
function is necessary and sufficient for uniform (asymptotic) stability in the large (on $\G$) of the origin.

Converse results for stability and uniform stability were initially developed by Krasovskii \cite{Kras55-PMM}
and Kurzweil \cite{Kurz55-CMJ}.  At the same time, without assuming uniqueness of solutions,
Yoshizawa presented a continuous Lyapunov function assuming stability of the origin in \cite{Yosh55-Kyoto}.
Smooth Lyapunov functions for stability and uniform stability of the origin, without the assumption of unique
trajectories, were provided by Kurzweil and Vrkoc \cite{KuVr57-CMJ}.

Almost simultaneously\footnote{Twice in the 1950's, similar results were submitted almost simultaneously.
Similar results were published by Krasovskii \cite{Kras55-PMM} (submitted 12 November 1954) and
Kurzweil \cite{Kurz55-CMJ} (submitted 2 December 1954).  Again involving Kurzweil, similar results
were published by Kurzweil \cite{Kurz56-AMST} (submitted 6 July 1955) and Massera \cite{Mass56-AM}
(submitted 30 August 1955).},
Kurzweil \cite{Kurz56-AMST} and Massera \cite{Mass56-AM} demonstrated that, when the righthand side
of \eqref{eq:sys} is continuous, if the origin is uniformly globally asymptotically stable\footnote{Kurzweil's result
is actually for uniform asymptotic stability in the large on the domain of attraction $\D \subseteq \R^n$.} then there exists
a smooth (infinitely differentiable) Lyapunov function.  While Massera assumed unique solutions to
\eqref{eq:sys}, Kurzweil did not\footnote{Despite the similar results on the existence of a smooth Lyapunov
function, Massera's assumption of unique solutions allowed a much shorter proof.  In fact, \cite{Mass56-AM}
grew out of a short course Massera provided in Varenna, Italy in 1954 and contains a nice survey of many topics in
stability theory.}.  In order to accommodate the lack of a unique solution, Kurzweil
extended Massera's construction \eqref{eq:Massera} by taking the supremum over all solutions
from the initial condition, $x \in \R^n$:
\begin{equation}
	\label{eq:Kurzweil}
	V(x,t) = \sup_{\phi \in \S(x)} \int_t^\infty \alpha(|\phi(\tau,t,x)|)  d\tau,
\end{equation}
where, as in \eqref{eq:Massera}, $\alpha: \R_{\geq 0} \rightarrow \R_{\geq 0}$ is an appropriately
chosen continuous function.  The functions defined by \eqref{eq:Massera} and \eqref{eq:Kurzweil}
can be shown to be locally Lipschitz, after which a transfinite smoothing procedure that maintains the
desired Lyapunov function properties, is applied.  As in \cite{Mass49-AM}, both \cite{Mass56-AM} and
\cite{Kurz56-AMST} show that if the righthand side of \eqref{eq:sys} is periodic in $t$ or independent
of $t$, then so is the derived Lyapunov function.  It is worth noting that for (uniform)
stability this does not hold; i.e., even for systems \eqref{eq:sys} that are independent of $t$ {\rev and that
possess a (uniformly) stable equilibrium point,} it is not
always possible to find a Lyapunov function that is independent of $t$ (see \cite[p.\ 46]{Kras63}).

By the end of the 1950's, answers to most of the converse questions for Theorem \ref{thm:ModernLyapunov},
including its in-the-large variants, had thus been obtained.  Subsequent research focused on more general
systems and on more general stability concepts.  {\rev The one remaining converse from
Theorem \ref{thm:ModernLyapunov} relates to equiasymptotic stability of the origin.  This converse appears
to have been originally derived by Hahn in \cite[Theorem 49.1]{Hahn67}, where, similar to the above
observation on (uniform) stability, even for systems independent of $t$ it is not always possible to find
a Lyapunov function that is independent of $t$.
}

\subsection{Extensions and consolidation - the 1960's}
Lyapunov's second method was extended to so-called ``general dynamical systems'', namely
dynamical systems axiomatically defined based on the attainability sets of differential equations
without unique solutions on metric spaces.  Research on such systems was initiated by Barbashin \cite{Barb48-Russian}
and converse Lyapunov theorems for systems with unique solutions were provided by Zubov \cite{Zubo64}.
The extension to systems without unique solutions was provided in a series of papers by Roxin
\cite{Roxi65a-JDE}, \cite{Roxi66-RCMP}, and \cite{Roxi66-SJC}, where a distinction was made between
the stability behavior of all solutions and the stability behavior of at least one solution.  As we describe
precisely in Definition \ref{def:KL_stable} below, these stability properties are termed {\em strong stability}
and {\em weak stability}, respectively.  An excellent summary of the initial work on general dynamical
systems can be found in \cite{Kloe78-MCT}.

As a particular case of both partial stability \cite{Voro05-AiT} and stability with respect to two
measures \cite{Movc60-PMM}, Hoppensteadt \cite{Hopp66-TAMS} derived a continuously differentiable Lyapunov function
for asymptotic stability of the origin for a parametrized non-autonomous differential equation where
the parameters take values in an unbounded set.  Wilson \cite{Wils69-TAMS} then extended this by deriving a smooth Lyapunov
function under the assumption of uniform asymptotic stability of a closed, but not necessarily bounded, set
and Lakshmikantham and Salvadori \cite{LaSa76-BUMI} provided a continuous Lyapunov function under the assumption of
stability with respect to two measures.

In the 1960's, five books appeared in English that summarized much of the available theory on Lyapunov's
second method, including results on the converse question \cite{Hahn63, Kras63, Zubo64, Yosh66, Hahn67}.
Along with the survey papers \cite{Anto58} and \cite{KaBe60-JBE} and the very readable
\cite{LaLe61}, these texts made Lyapunov's methods widely accessible
to the West\footnote{{\rev It is tempting to speculate that the wealth of translated material in the 1960's
was directly due to the initial accomplishments of the Soviet Union in the space race with the
launch of the Sputnik satellite in October 1957 and Yuri Gagarin's orbital flight in April 1961.
In particular, what appear to be hasty translations by government departments were made
of the 1951 text of Lur'e \cite{Lure51} in 1957 in the United Kingdom by Her Majesty's Stationery Office,
and of the 1956 text of
Malkin \cite{Malk56} in 1959 by the United States Atomic Energy Commission.  Other important
Soviet texts that were translated around this time include the 1955 text of Letov \cite{Leto61}
(translated 1961), the 1956 text of Chetaev \cite{Chet61} (translated 1961), the 1957 text
of Zubov \cite{Zubo64} (translated 1964), and the 1959 text of Krasovskii \cite{Kras63} (translated 1963).
Furthermore, the preeminent Russian language control journal published since 1936,
{\cyr Prikladnaya Matematika i Mekhanika} (Prikladnaya Matematika i Mekhanika), was regularly translated as the {\em Journal of
Applied Mathematics and Mechanics} beginning in 1958.  While Cold War tensions, and the
space race in particular, is likely one driver behind this rush of translations, a level of translation activity that
has not been seen since,
it is also worth noting that the International Federation of Automatic Control was formed in 1957
with a goal of international scientific exchange and whose first two presidents came from
the United States of America (H.~Chestnut, 1957--1959) and the Soviet Union
(A.~L.~Letov, 1959--1961).  Historical narratives are rarely simple.}
}.  Of particular note, similar to the candidate Lyapunov function he used in \cite{Yosh55-Kyoto},
Yoshizawa proved his converse theorems in \cite{Yosh66} based on a function defined as
\begin{equation}
	\label{eq:Yoshizawa}
	V(x,t) = \sup_{\phi \in \S(x), \ \tau \geq t} \alpha(|\phi(\tau,t,x)|)e^{c\tau}
\end{equation}
where $\alpha \in \Kinf$ and $c \in R_{>0}$ are chosen appropriately.

Krasovskii \cite{Kras63}, on the other hand, used a different proof technique to that initially
proposed by Massera in \cite{Mass49-AM}.  Of particular note is that Krasovskii's technique
allowed him to derive a converse theorem not just for asymptotic stability, but also
for Lyapunov's first instability theorem (see Section \ref{sec:Instability}).  The core of this technique
rests on what Krasovskii labels ``Property A'' \cite[Definition 4.1]{Kras63}:

\begin{quote}
{\it Property A:} Let $\{h_k\}_{k=0}^\infty$ be a monotonically decreasing sequence satisfying
\begin{equation}
	h_0 = |0|_{\G}, \qquad
	\lim_{k \rightarrow \infty} h_k = 0 .
\end{equation}
Suppose that for every closed bounded region $H \subset \R^n$ satisfying $0 \in \overline{H} \subset \G$,
and for every $k > 0$
there is a number $T_k$ such that whenever $t_0 \geq T_k$, $x \in \G \backslash \B_{h_k}$,
there exists $t \in [t_0 - T_k, t_0 + T_k]$ so that $\phi(t, t_0, x) \notin H \backslash \{0\}$.
\end{quote}

In words, for any neighborhood of the origin in $\G$ and for a sequence of decreasing balls (centered at
the origin), solutions cannot stay in $H$ indefinitely.  Intuitively, this property is satisfied for both
asymptotically stable and unstable equilibria -- though not necessarily for stable equilibria.
Krasovskii then demonstrated \cite[Theorem 4.3]{Kras63} that Property A is equivalent to the existence
of a function $V$ that is decrescent and such that its derivative along solutions of \eqref{eq:sys} is
sign-definite.  This provides one of the main technical tools in \cite{Kras63} to derive converse theorems
for asymptotic stability and for instability (see Section \ref{sec:Instability} below for the latter).  An interesting
result, apparently not available elsewhere, is a converse theorem for equiasymptotic stability of the origin
\cite[Theorem 10.2]{Kras63} where
the total derivative of the Lyapunov function is negative semidefinite and the supremum of the total
derivative integrates to negative infinity.

\subsection{1970's onward}

In the late 1970s researchers began to examine differential and difference inclusions.
Roxin \cite{Roxi65b-JDE} demonstrated how differential
inclusions (also called contingent equations) give rise to the general dynamical systems of Barbashin
\cite{Barb48-Russian}, and so the specific results for differential inclusions could be seen as a special
case of Roxin's results.  However, the strength of Lyapunov's second method has generally been that
one need not generate system trajectories whereas the general systems approach of Roxin requires
knowledge of the attainability function or the set of solutions.  The specialization to difference
and differential inclusions allows the formulation of decrease conditions that depend on the
set-valued mapping defining the inclusion rather than requiring knowledge of solutions.

Meilakhs \cite{Meil78-AiT} derived a continuously differentiable Lyapunov function given uniform asymptotic
stability of the origin for all solutions of a differential inclusion derived from a parametrized differential equation where the
parameters vary over a closed bounded linearly connected set.  Molchanov and Pyatnitskii studied the
problem of absolute stability {\rev described in Section~\ref{sec:Lure}}.  In
\cite{MoPy86a-AiT} and \cite{MoPy86b-AiT} they formulated the Lur'e problem as a stability problem for a
differential inclusion and demonstrated the existence of a Lyapunov function of an approximate quadratic
form.  Similar to the algebraic criteria of Section \ref{sec:constructive}, in \cite{MoPy86c-AiT} and \cite{MoPy89-SCL}
Molchanov and Pyatnitskii then derived necessary and sufficient criteria for a Lyapunov function in terms of solvability of certain
matrix equations.

A result on the existence
of a so-called {\em control Lyapunov function} under the assumption of asymptotic controllability to the origin was
provided by Sontag \cite{Sont83-SICON}.  This is closely related to the converse question for weak asymptotic stability
of the origin for differential inclusions, which Smirnov answered in \cite{Smir90a-AiT} and
\cite{Smir90b-AiT} for differential inclusions described by convex processes.    Converse theorems for both weak and strong
stability of time-varying differential inclusions defined on a real Banach space were provided by
Deimling \cite[Propositions 14.1 and 14.2]{Deim92}.

A converse theorem for uniform global asymptotic stability of a compact set for a differential inclusion
under fairly weak assumptions was provided by Lin et al.\ \cite{LSW96-SICON} and converse theorems
for both uniform global strong and weak asymptotic stability were provided by Clarke et al.\ \cite{CLS98-JDE}
where a particular impediment in the weak case was identified.  We discuss this impediment
in Section \ref{sec:WeakKLStable}.

The first converse Lyapunov theorems for discrete time systems described by non-autonomous ordinary difference
equations were derived by Gordon \cite{Gord72-SJC} for stability and uniform asymptotic stability of the origin.
A converse theorem for strong uniform global asymptotic stability
of a difference inclusion was provided by Jiang and Wang \cite{JiWa02-SCL}.

{\rev Similar to previous results that time-independent or periodic systems yield time-independent or
periodic Lyapunov functions, Rosier \cite{Rosi92-SCL} demonstrated that for homogeneous systems with non-unique
solutions, an asymptotically stable origin implies the existence of a smooth homogeneous Lyapunov function.
}

Converse theorems for both strong and weak stability with respect to two measures for
both differential and difference inclusions were provided by Teel and Praly \cite{TePr00-CoCV}
and by the author and Teel in
\cite{KeTe04b-SCL}, \cite{KeTe04-SICON}, \cite{KeTe05-SICON}, and \cite{KeTe07-MCSS}.
As these general results subsume
much previous work as special cases, we specifically survey these results in Section \ref{sec:KellettTeel}.

Recently, Karafyllis and Tsinias \cite{KaTs03-SICON}  and Karafyllis \cite{Kara05-IMA}
developed converse theorems
for strong equiasymptotic stability of the origin for differential and difference inclusions arising from
perturbed difference and differential equations.  Rather than equiasymptotic stability they use the
terminology non-uniform in time stability.

{\rev Also recently, Kloeden  and co-authors \cite{Kloe98,Kloe00,GKSW07-DCDS} noted that nonautonomous
systems naturally give rise to nonautonomous invariant sets.  This then leads to three notions of attractor
and stability, refered to as pullback, forward, and uniform attractor/stability.  Similar to Massera's result
\cite[Theorem 7]{Mass56-AM} for systems periodic in, or independent of, time the definitions of pullback, forward,
and uniform attractor coincide.  Appropriate Lyapunov functions were defined and converse results
presented in \cite[Theorem 29]{GKSW07-DCDS} for pullback, forward, and uniform attractors/stability of
nonautonomous differential equations, while \cite{Kloe00} presents a converse result for pullback attraction of
a nonautonomous difference equation.  The constructions used in these references are similar to that proposed
by Yoshizawa \eqref{eq:Yoshizawa}. }

% -----------------------------------------------------------------------------------------------------------------------------------------
%		 KL-STABILITY WITH RESPECT TO TWO MEASURES
% -----------------------------------------------------------------------------------------------------------------------------------------

\section{$\KL$-stability with respect to two measures for difference and differential inclusions}
\label{sec:KellettTeel}
For a set-valued map $F(\cdot)$ we use the notation $F: \R^n \rightrightarrows \R^n$ to denote
that $F(\cdot)$ maps points in $\R^n$ to subsets of $\R^n$.  In this section, for comparative purposes,
we present some specific converse theorems for difference inclusions
\begin{equation}
	\label{eq:DT_Incl}
	x^+ \in F(x)
\end{equation}
and differential inclusions
\begin{equation}
	\label{eq:CT_Incl}
	\dot{x} \in F(x)
\end{equation}
where  $F: \G \rightrightarrows \G$ for \eqref{eq:DT_Incl} and $F: \G \rightrightarrows \R^n$ for \eqref{eq:CT_Incl}, and
where $\G \subset \R^n$.
In an abuse of notation, in order to avoid unnecessary duplication in the results that follow, we use $t$ both as
$t \in \R_{\geq 0}$ when refering to the continuous time system \eqref{eq:CT_Incl} and
$t \in \Z_{\geq 0}$ when refering to the discrete time system \eqref{eq:DT_Incl}.

For completeness, we provide here regularity definitions for set-valued maps.
{\rev Note that, for sets $A,B \subset \R^n$, $A+B \subset \R^n$ denotes the Minkowski sum.}
\begin{definition}
	Let $\mathcal{O} \subset \R^n$ be open.  The set-valued map $F: \R^n \rightrightarrows \R^n$ is:
	\begin{itemize}
		\item {\em upper semicontinuous on $\mathcal{O}$} if for each $x \in \mathcal{O}$ and $\varepsilon > 0$
			there exists $\delta > 0$ such that for all $\xi \in \mathcal{O}$ satisfying $|x-\xi| < \delta$ we have
			$F(\xi) \subset F(x) + \B_\varepsilon$;
		\item {\em continuous on $\mathcal{O}$} if, in addition to being upper semicontinuous on $\mathcal{O}$, for each
			$x \in \mathcal{O}$ and $\varepsilon > 0$ there exists $\delta > 0$ such that, for
			$\xi \in \mathcal{O}$ satisfying $|\xi-x| < \delta$ we have $F(x) \subset F(\xi) + \B_\varepsilon$; and
		\item {\em locally Lipschitz on $\mathcal{O}$} if for each $x \in \mathcal{O}$ there exists a
			neighborhood $\mathcal{U} \subset \mathcal{O}$ of $x$ and $L > 0$ such that
			$x_1, x_2 \in \mathcal{O}$ implies $F(x_1) \subset F(x_2) + L|x_1 - x_2| \overline{\B}$.
	\end{itemize}
\end{definition}

Note that the concept of upper semicontinuity for a set-valued map is not the same as that for extended real-valued functions. In fact, for $f:\R^n \rightarrow \R^n$, the set-valued map $x \mapsto \{f(x)\}$ is upper semicontinuous if and only if the
extended real-valued function $x \mapsto f(x)$ is continuous.

The results of this section generally require a common set of assumptions with regards to the set-valued
map defining the difference or differential inclusion.
\begin{definition}
	The set-valued map $F:\G \rightrightarrows \G$ satisfies the {\em discrete time basic conditions}
	if, on $\G$, it has nonempty and compact values, and is upper semicontinuous.
\end{definition}
\begin{definition}
	The set-valued map $F:\G \rightrightarrows \R^n$ satisfies the {\em continuous time basic conditions}
	if, on $\G$, it has nonempty, compact, and convex values, and is upper semicontinuous.
\end{definition}
The continuous time basic conditions are essentially required in order to guarantee existence of solutions
to the differential inclusion (see \cite{Fili88}).  These conditions also provide certain technical properties on
the solution sets.  By contrast, solutions to the difference inclusion \eqref{eq:DT_Incl} will exist so long as the
mapping is nonempty.  However, the discrete time basic conditions enable certain technical results such as
closeness of solutions properties (see \cite{KeTe05-SICON}).

For systems that do not give rise to unique solutions there are two natural stability
notions that were identified by Roxin \cite{Roxi65a-JDE}.  The first is the property that {\em all} solutions
must satisfy a desired stability estimate while the second is the property that {\em at least one} solution
must satisfy a desired stability estimate.  Roxin termed these properties ``strong stability'' and ``weak stability'',
respectively.

In the framework of difference and differential inclusions, the subsequent results subsume many
commonly encountered system models including ordinary difference and differential equations, such systems
with discontinuous righthand sides, and controlled or perturbed systems.  To further extend the reach of these
results, we can consider a generalization of uniform global asymptotic stability that was introduced by
Movchan \cite{Movc60-PMM} refered to as stability with respect to two measures or stability with respect
to two metrics.
\begin{definition}
	\label{def:KL_stable}
	Let $\omega_i : \G \rightarrow \R_{\geq 0}$, $i=1,2$, be continuous functions.
	We say that \eqref{eq:DT_Incl} (or \eqref{eq:CT_Incl}) is strongly $\KL$-stable with respect to
	$(\omega_1, \omega_2)$ if {\rev (\eqref{eq:CT_Incl} is forward complete and)} there exists $\beta \in \KL$ such that for every initial condition
	$x \in \G$
	all solutions $\phi \in \S(x)$ satisfy
	\begin{equation}
		\label{eq:2meas_stable}
		\omega_1(\phi(t,x)) \leq \beta(\omega_2(x),t), \quad \forall t \in \Z_{\geq 0} \ (\forall t \in \R_{\geq 0}).
	\end{equation}
	We say that \eqref{eq:DT_Incl} (or \eqref{eq:CT_Incl}) is weakly $\KL$-stable with respect to
	$(\omega_1, \omega_2)$ if the above property holds for at least one solution $\phi \in \S(x)$.
\end{definition}

Note that, in some sense, both ``stability with respect to two measures'' and ``stability with respect to
two metrics'' are unsatisfactory terminology as the functions $\omega_i$ are neither measures nor metrics in the
usual mathematical sense of measure or metric.  Nonetheless, the usage has become standard and
we will use the terminology ``stability with respect to two measures''.

{\rev Observe that forward completeness is only explicitly required for continuous time strong $\KL$-stability
with respect to two measures as this is not guaranteed {\em a priori} by the stability estimate.  Forward completeness
of difference inclusions is guaranteed by virtue of the set-valued mapping taking points in $\G$ to subsets of
$\G$.  Finally, in the case of continuous time weak $\KL$-stability, since we are not necessarily interested in the
behavior of all solutions, it may in fact be the case that some solutions
cannot be continued for all time.
}

$\KL$-stability with respect to two measures is a generalization of uniform asymptotic stability of the
origin in the large on $\G$.  This is the case where $\omega_1(x) = \omega_2(x) = |x|$ and, hence, when
additionally $\G = \R^n$, we see that $\KL$-stability with respect to $(|\cdot|, |\cdot|)$ is, in fact, uniform
global asymptotic stability of the origin.  Additionally, this stability property also encompasses uniform
asymptotic stability of a closed set $\A \subset \G$ by taking $\omega_1(x) = \omega_2(x) = |x|_{\A}$
where $|x|_{\A} = \min_{y \in \A} |x-y|$.  Many other examples, such as output stability and stability of
a particular trajectory are possible.

As a particular example, it is possible to deal with non-autonomous systems \eqref{eq:sys} as autonomous
systems by the technique of state augmentation; that is, consider states $x = (z,t) \in \R^n \times \R_{\geq 0}$
and
\begin{equation}
	\label{eq:augmented}
	\begin{array}{rl}
		\dot{z} & = f(z,t) \\
		\dot{t} & = 1 .
	\end{array}
\end{equation}
Then uniform global asymptotic stability of the origin is equivalent to $\KL$-stability with respect to
$(\omega, \omega)$ where $\omega(z,t) = |z|$ for all $(z,t) \in \R^n \times \R_{\geq 0}$.
A similar approach can be taken for systems \eqref{eq:sys}
parametrized by a parameter vector, $\theta \in \R^m$, where the system equations can be augmented
by $\dot{\theta} = 0$.  Note that a limitation of this technique is that it is necessary to impose 
regularity conditions on the $t$ or $\theta$-dependence of $f$ than are strictly required.

As another example, consider the second order system
\begin{equation}
	\label{eq:example}
	\begin{array}{rcl}
		\dot{x}_1 & = & x_2 + (1 - x_1^2 - x_2^2)x_1 \\
		\dot{x}_2 & = & -x_1 + (1 - x_1^2 - x_2^2)x_2.
	\end{array}
\end{equation}
This system has the origin as an unstable equilibrium point and the unit circle as a
uniformly asymptotically stable periodic orbit.  For $x \in \G = \R^2 \backslash \{0\}$, define the function
\begin{equation}
	\label{eq:example_meas}
	\omega(x) := \left\{\begin{array}{cl}
		\frac{1 - |x|}{|x|}, & x \in \B \backslash \{0\} \\
		|x| - 1, & x \in \R^2 \backslash \B.
	\end{array} \right.
\end{equation}
Then \eqref{eq:example} is $\KL$-stable with respect to $(\omega, \omega)$ which captures the system
behavior both in terms of the unstable equilibrium as well as the asymptotically stable periodic orbit.

\subsection{Converse theorems for strong $\KL$-stability}
In the context of differential equations with unique solutions,
Massera observed that certain stability properties, namely equiasymptotic stability and uniform
asymptotic stability, have an inherent robustness property
in that the set of solutions is an open set \cite[Theorem 8]{Mass56-AM}.  In other words, near any solution that
satisfies an equiasymptotic or uniform asymptotic stability estimate, there are other solutions that also satisfy that
estimate.  In the case of strong $\KL$-stability
with respect to two measures for both difference and differential inclusions, as a first
step towards various converse theorems, we make a connection between
robust stability and the existence of a smooth Lyapunov function.  Then, to complete
a converse Lyapunov theorem, we present various conditions that guarantee robust stability.

For both difference and differential inclusions we define robust stability in terms of stability of
a perturbed inclusion.  Define
\[ \A := \left\{ x \in \G : \sup_{t \in \mathbb{T}, \phi \in \S(x)} \omega_1(\phi(t,x)) = 0 \right\} \]
where $\mathbb{T} = \Z_{\geq 0}$ for \eqref{eq:DT_Incl} and $\mathbb{T} = \R_{\geq 0}$
for \eqref{eq:CT_Incl}.
For continuous functions $\sigma, \delta: \G \rightarrow \R_{\geq 0}$, such that
 $\sigma(x), \delta(x) > 0$, for all $x \in \G \backslash \A$, and
\[ \{x\} + \sigma(x)\overline{\B} \subset \G, \quad \{x\} + \delta(x)\overline{\B} \subset \G \]
we define the perturbed inclusions
\begin{align}
	x^+ \in F_\sigma (x) & := \left\{ v \in \R^n : v \in \{ \eta \} + \sigma(\eta) \overline{\B}_n, \ \eta \in
		F\left(x + \sigma(x) \overline{\B}_n\right) \right\}, \quad {\rm and}	\label{eq:DT_perturb} \\
	\dot{x} \in F_\delta (x) & := \overline{\rm co} F\left(\{x\} + \delta(x) \overline{\B}\right) + \delta(x) \overline{\B}.
		\label{eq:CT_perturb}
\end{align}
Note that since differential inclusions deal with infinitesimals it is possible to define the inner and outer
perturbations of \eqref{eq:CT_perturb}
on the basis of the same point $x \in \G$.  This is in contrast to the perturbed difference inclusion \eqref{eq:DT_perturb}
where the outer perturbation needs to be a superset of the set-valued map applied to the inner perturbation.
Also note that it is necessary to take the closed convex hull in \eqref{eq:CT_perturb} to ensure that
$F_\delta$ satisfies the continuous time basic conditions.

If \eqref{eq:DT_perturb}, respectively \eqref{eq:CT_perturb}, is $\KL$-stable with respect to
$(\omega_1, \omega_2)$ then we say that \eqref{eq:DT_Incl}, respectively \eqref{eq:CT_Incl},
is robustly $\KL$-stable with respect to $(\omega_1, \omega_2)$.

\begin{theorem}{\cite[Theorem 2.7]{KeTe05-SICON}}
	\label{thm:KeTe_main}
	Let $F: \G \rightrightarrows \G$ satisfy the discrete time basic conditions on $\G$.
	The difference inclusion \eqref{eq:DT_Incl} is robustly $\KL$-stable with respect to $(\omega_1, \omega_2)$ if
	and only if there exists a smooth Lyapunov function with respect to $(\omega_1, \omega_2)$ on $\G$; i.e.,
	a smooth function $V: \G \rightarrow \R_{\geq 0}$ and $\alpha_1, \alpha_2 \in \Kinf$
	such that for all $x \in \G$
	\begin{equation}
		\alpha_1(\omega_1(x)) \leq V(x) \leq \alpha_2(\omega_2(x))
	\end{equation}
	\begin{equation}
		\label{eq:DT_decrease}
		\max_{f \in F(x)} V(f) \leq V(x)e^{-1}
	\end{equation}
	where $e^{-1}$ is the exponential function evaluated at $-1$.
\end{theorem}

\begin{theorem}{\cite[Theorem 1]{TePr00-CoCV}}
	\label{thm:TePr_main}
	Let $F: \G \rightrightarrows \R^n$ satisfy the basic conditions on $\G$.
	The differential inclusion \eqref{eq:CT_Incl} is robustly $\KL$-stable with respect to $(\omega_1, \omega_2)$
	if and only if \eqref{eq:CT_Incl} is
	forward complete on $\G$ and there exists a smooth Lyapunov function with respect to $(\omega_1, \omega_2)$
	on $\G$; i.e., a smooth function $V: \G \rightarrow \R_{\geq 0}$ and $\alpha_1, \alpha_2 \in \Kinf$
	such that for all $x \in \G$
	\begin{equation}
		\label{eq:UpperLower}
		\alpha_1(\omega_1(x)) \leq V(x) \leq \alpha_2(\omega_2(x))
	\end{equation}
	\begin{equation}
		\label{eq:CT_decrease}
		\max_{w \in F(x)} \left\langle \tfrac{\partial}{\partial x} V(x), w \right\rangle \leq - V(x).
	\end{equation}
\end{theorem}

Observe that forward completeness is explicitly required for differential inclusions whereas this is
guaranteed for difference inclusions by virtue of the fact that the set-valued map takes points in
$\G$ to subsets of $\G$; i.e., by definition, solutions to \eqref{eq:DT_Incl} cannot escape $\G$.
{\rev Demonstrating forward completeness for differential inclusions can be accomplished via Lyapunov
methods \cite{AnSo99-SCL}.}

The decrease conditions \eqref{eq:DT_decrease} and \eqref{eq:CT_decrease} guarantee that
the Lyapunov functions decrease exponentially along solutions of \eqref{eq:DT_Incl} and
\eqref{eq:CT_Incl}, respectively, where the constant $e^{-1}$ is chosen to mirror the exponential
decrease implied by \eqref{eq:CT_decrease}.  Given any Lyapunov function that does not decrease exponentially,
it is always possible to find a nonlinear scaling such that the nonlinear scaling of the Lyapunov function
is also a Lyapunov function and decreases exponentially.  It is important to note that this does not imply that
system solutions decrease exponentially.  In fact, a sufficient condition for the exponential
decrease of solutions (i.e., exponential stability) is that the decrease condition \eqref{eq:CT_decrease} is
satisfied and that the upper and lower bounds \eqref{eq:UpperLower}
be quadratic.

In order to compare the Lyapunov function given by \eqref{eq:UpperLower}-\eqref{eq:CT_decrease} and
that in Theorem \ref{thm:ModernLyapunov} we briefly examine \eqref{eq:augmented}.  For uniform global
asymptotic stability of the origin we observed that, for all $x = (z,t) \in \R^n \times \R_{\geq 0}$,
 $\omega_1(x) = \omega_2(x) = |z|$.  Therefore, for \eqref{eq:augmented}, \eqref{eq:UpperLower} is
\[ \alpha_1(|z|) \leq V(z,t) \leq \alpha_2(|z|) \]
which is precisely items (\ref{it:pos_def}) and (\ref{it:decrescent}) of Theorem \ref{thm:ModernLyapunov}.
Furthermore, \eqref{eq:CT_decrease} is
\[ \max_{w = [f(z,t)^T \ 1]^T} \left\langle \left[ \tfrac{\partial V}{\partial z}^T \ \tfrac{\partial V}{\partial t}  \right]^T, w \right\rangle
	=  \tfrac{\partial V}{\partial z}^T f(z,t)  + \tfrac{\partial V}{\partial t} \leq - V(x) \leq - \alpha_1(|z|) \]
which is precisely item (\ref{it:decrease}) of Theorem \ref{thm:ModernLyapunov} where
$\rho = \alpha_1 \in \Kinf$.  Therefore, item (\ref{it:UGAS}) of Theorem \ref{thm:ModernLyapunov} implies
that $z=0$ is uniformly globally asymptotically stable for \eqref{eq:augmented}.

It remains an open question for both difference and differential inclusions as to whether or not $\KL$-stability with
respect to $(\omega_1, \omega_2)$ is generally robust.  However, sufficient conditions for robustness
have been demonstrated in several special cases.  Many of these conditions are similar between discrete and
continuous time.  One such condition is related to the regularity of the set-valued map.
\begin{theorem}{\cite[Theorem 2.10]{KeTe05-SICON}}
	\label{thm:KeTe_cts}
	Let $F: \G \rightrightarrows \G$ satisfy the discrete time basic conditions on $\G$ and be continuous
	on an open set containing $\G \backslash \A$.  If \eqref{eq:DT_Incl} is strongly $\KL$-stable with respect to
	$(\omega_1, \omega_2)$ on $\G$, then it is robustly $\KL$-stable with respect to $(\omega_1, \omega_2)$
	on $\G$.
\end{theorem}
\begin{theorem}{\cite[Theorem 2]{TePr00-CoCV}}
	\label{thm:TePr_cts}
	Let $F: \G \rightrightarrows \R^n$ satisfy the continuous time basic conditions on $\G$ and be locally Lipschitz
	on an open set containing $\G \backslash \A$.  If \eqref{eq:CT_Incl} is strongly $\KL$-stable with respect to
	$(\omega_1, \omega_2)$ on $\G$, then it is robustly $\KL$-stable with respect to $(\omega_1, \omega_2)$
	on $\G$.
\end{theorem}

Note that Theorems \ref{thm:TePr_main} and \ref{thm:TePr_cts}
imply the existence of a Lyapunov function on $\R^2 \backslash \{0\}$ for \eqref{eq:example} with respect to the
measurement function of \eqref{eq:example_meas}.  Such a Lyapunov function is closely related to Conley's
complete Lyapunov functions \cite{Conl78} which are defined on the entire space.

Other sufficient conditions for robust $\KL$-stability of difference and differential inclusions
have been provided in \cite{KeTe05-SICON}, \cite{KeTe07-MCSS}, and \cite{TePr00-CoCV}.
In particular, for a compact attractor $\A \subset \R^n$, when $\omega_1 = \omega_2 = \omega$ is a
proper indicator function\footnote{A proper indicator for $\A$ on $\D$ is a continuous function
$\omega: \R^n \rightarrow \R_{\geq 0}$ such that $\omega(x) = 0$ for $x \in \A$, $\omega(x) > 0$ for $x \in \D \backslash \A$,
and $\lim_{x \rightarrow \partial \D} \omega(x) = c$ for some $c \in \R_{>0} \cup \infty$.} for $\A$ on its domain of attraction $\D \subset \R^n$, the basic conditions are sufficient to guarantee robustness; i.e., the extra regularity of Theorems
\ref{thm:KeTe_cts} and \ref{thm:TePr_cts} is not required.  An additional sufficient condition for robustness in the case
of a single measurement function is related to how solutions behave in reverse time; see
\cite[Theorem 4]{KeTe07-MCSS} and \cite[Theorem 3]{TePr00-CoCV}.  Finally,  for difference
inclusions where the set-valued mapping has no specific regularity requirement but is compact and nonempty,
the existence of a continuous
Lyapunov function as described in Theorem \ref{thm:KeTe_main} is sufficient to guarantee robust
$\KL$-stability of \eqref{eq:DT_Incl} (see \cite[Theorem 2.8]{KeTe05-SICON}).

In concluding our discussion on converse theorems for strong $\KL$-stability with respect to two measures,
we mention that results similar to those described above are available in the framework of hybrid systems that are defined by both
difference and differential inclusions and particular rules about when solutions evolve according to
the difference inclusion \eqref{eq:DT_Incl} or the differential inclusion \eqref{eq:CT_Incl}; see \cite{CTG07-TAC},
\cite{CTG08-TAC}, and \cite[Section 7.5]{GST12}.

\subsection{Converse theorems for weak $\KL$-stability}
\label{sec:WeakKLStable}
As described above, there are many cases in which strong $\KL$-stability with respect to
two measures is in fact robust and, consequently, several converse Lyapunov theorems are possible.
Fewer results have been obtained in the case of weak $\KL$-stability with respect to two measures.
In fact, currently available results are limited to the case where both measurement functions
are given by the distance to a closed (possibly unbounded) set $\A \subset \R^n$; i.e.,
$\omega_1(x) = \omega_2(x) = |x|_{\A}$.  In part, this is
due to the fact that the currently available proofs depend on the measures being the same and
on the measures satisfying the triangle inequality.

\begin{theorem}{\cite[Theorem 6]{KeTe04b-SCL}}
	\label{thm:weak_dt}
	Suppose $F:\R^n \rightrightarrows \R^n$ satisfies the discrete time basic conditions, is
	continuous on $\R^n \backslash \A$, and that \eqref{eq:DT_Incl} is weakly $\KL$-stable
	with respect to $(|\cdot |_{\A}, |\cdot |_{\A})$.  Then there exists a weak discrete time
	Lyapunov function; that is, a smooth function $V:\R^n \rightarrow \R_{\geq 0}$ and
	$\alpha_1, \alpha_2 \in \Kinf$ such that, for all $x \in \R^n$
	\begin{equation}
		\alpha_1(|x|_{\A}) \leq V(x) \leq \alpha_2(|x|_{\A}), \quad {\rm and}
	\end{equation}
	\begin{equation}
		\min_{f \in F(x)} V(f) \leq V(x)e^{-1} .
	\end{equation}
\end{theorem}

The following theorem first appeared in \cite{KeTe00-MTNS} {\rev and requires the following assumption:

\begin{assumption}
	\label{as:1}
	For each $r \in \R_{>0}$ there exists $M_r \in \R_{>0}$ such that $|x|_\A \leq r$ implies
	$\sup_{w \in F(x)} |w| \leq M_r$.
\end{assumption}
}

\begin{theorem}{\cite[Theorem 2.1]{KeTe04-SICON}}
	\label{thm:weak_ct}
	Suppose $F:\R^n \rightrightarrows \R^n$ satisfies the continuous time basic conditions, is
	locally Lipschitz on $\R^n \backslash \A$, {\rev satisfies Assumption~\ref{as:1},} and that \eqref{eq:CT_Incl} is weakly $\KL$-stable
	with respect to $(|\cdot |_{\A}, |\cdot |_{\A})$.  Then there exists a weak continuous time
	Lyapunov function; that is, a locally Lipschitz function $V: \R^n \rightarrow \R_{\geq 0}$ and
	$\alpha_1, \alpha_2 \in \Kinf$ such that, for all $x \in \R^n$
	\begin{equation}
		\alpha_1(|x|_{\A}) \leq V(x) \leq \alpha_2(|x|_{\A}), \quad {\rm and}
	\end{equation}
	\begin{equation}
		\min_{w \in F(x)} DV(x;w) \leq - V(x),
	\end{equation}
	where $DV(x;w)$ denotes the Dini derivative at $x \in \R^n$ in the direction $w \in F(x)$.
\end{theorem}

Independent to \cite{KeTe00-MTNS} and \cite{KeTe04-SICON}, using techniques from optimal control,
Rifford \cite{Riff00-SICON}, under similar assumptions to those in Theorem \ref{thm:weak_ct} plus a linear growth condition
on the set-valued map,
derived a locally Lipschitz weak continuous time Lyapunov function for weak uniform asymptotic stability
of a compact set.

There is an interesting contrast between the continuous time case of Theorem \ref{thm:weak_ct}
where a locally Lipschitz Lyapunov function is obtained, versus the other presented cases
for strong stability in continuous time and both weak and strong stability in discrete time, where
smooth Lyapunov functions are obtained.  In fact, Clarke et al.\ \cite{CLS98-JDE} demonstrated that,
in general, it is not possible to find a smooth weak continuous time Lyapunov function for weak asymptotic stability of
\eqref{eq:CT_Incl}.  However, they provided a necessary condition in the form of a restriction on the
set-valued map satisfying a covering condition near the origin for such a smooth Lyapunov function.

{\rev
\begin{theorem}{\cite[Theorem 6.1]{CLS98-JDE}}
	\label{thm:covering}
	Suppose $F:\R^n \rightrightarrows \R^n$ satisfies the continuous time basic conditions and
	there exists a continuously differentiable weak Lyapunov function; i.e., a continuously differentiable
	function $V:\R^n \rightarrow \R_{\geq 0}$ and functions $\alpha_1, \alpha_2 \in \Kinf$ such that, for all
	$x \in \R^n$
	\begin{equation}
		\alpha_1(|x|) \leq V(x) \leq \alpha_2(|x|), \quad {\rm and}
	\end{equation}
	\begin{equation}
		\min_{w \in F(x)} \left\langle \tfrac{\partial}{\partial x} V(x), w \right\rangle \leq - V(x).
	\end{equation}
	Then, for any $\gamma \in \R_{>0}$ there exists $\Delta \in \R_{>0}$ such that
	\begin{equation}
		\B_\Delta \subset F\left(\B_\gamma\right) := \cup_{x \in \B_\gamma} F(x).
	\end{equation}
\end{theorem}

A simple example demonstrates that, in the case of continuous time weak $\KL$-stability, even if the above covering condition
is satisfied, a continuously differentiable weak Lyapunov function may fail to exist. Consider a system defined on $\R^2$ by
\begin{equation}
	\label{eq:2d_ex}
	\dot{x} \in \overline\B, \quad x \in \R^2.
\end{equation}
It is straightforward to see that $\A := \left\{(x_1, x_2) \in \R^2 : x_1^2 + x_2^2 = 1\right\}$ is
weakly $\KL$-stable with respect to $(| \cdot |_\A, | \cdot |_\A)$ for \eqref{eq:2d_ex}.  Suppose
it were possible to find a continuously differentiable function $V:\R^n \rightarrow \R_{\geq 0}$
and $\alpha_1, \alpha_2 \in \Kinf$ so that
	\begin{equation}
		\label{eq:2d_bnds}
		\alpha_1(|x|_\A) \leq V(x) \leq \alpha_2(|x|_\A), \quad {\rm and}
	\end{equation}
	\begin{equation}
		\min_{w \in \overline\B} \left\langle \tfrac{\partial}{\partial x} V(x), w \right\rangle \leq - V(x).
	\end{equation}
The above implies that $\tfrac{\partial}{\partial x} V(x) \neq 0$ for all $x \in \R^2 \backslash \A$.  On the other hand,
since $V$ is continuously differentiable, it obtains its minimum and maximum on $\overline\B$.
Equation \eqref{eq:2d_bnds} implies that $V$ obtains its minimum everywhere on the boundary
of $\overline\B$.  Consequently, $V$ must obtain its maximum on the interior of $\overline\B$,
contradicting that $\tfrac{\partial}{\partial x} V(x) \neq 0$ for all $x \in \R^2 \backslash \A$.  Hence, a continuously
differentiable weak Lyapunov function cannot exist for \eqref{eq:2d_ex}.
}

The covering condition of Theorem~\ref{thm:covering} is related to a similar covering condition derived by Brockett \cite{Broc82-BCh}
in the context of designing continuous
feedback stabilizers for controlled differential equations.  An example of a system that does not
satisfy this covering condition is a tricycle that needs to be steered to the origin where there are
clearly initial configurations that require a discontinuous decision to be made in terms of turning
the handlebars left or right.  This example belongs to the general class of systems refered to as
nonholonomic systems.

Such discontinuous feedback stabilizers suffer from a lack of robustness.  In particular, for the initial
configurations of the tricycle where a discontinuous decision must be made, arbitrarily small errors
in measuring the configuration can lead to a so-called chattering phenomenon near the point of
discontinuity and the system therefore never approaches the origin (see \cite{Sont99-COCV} for
an extended discussion).  Given the connections between smooth Lyapunov functions and
robustness, it is not surprising then that a smooth Lyapunov function is not possible for weak asymptotic
stability of \eqref{eq:CT_Incl}.

However, having made that observation, it is perhaps surprising that a smooth weak discrete time Lyapunov function is
possible for weak asymptotic stability of \eqref{eq:DT_Incl}.  This stems from the fact that a sampled-data
or discrete time implementation can circumvent the lack of robustness just described.  In the tricycle
example, a decision is made near the point of discontinuity and that decision is adhered to for a set
period of time.  This moves the system far enough away from the point of discontinuity that the aforementioned
chattering phenomenon does not occur, where there is clearly a relationship between how large the
measurement errors are and how long the set period of time is.  In some generality, the construction
of robust discontinuous feedback stabilizers was dealt with in
\cite{CLRS00-SICON}, \cite{CLSS97-TAC}, \cite{KeTe00-CDC}, \cite{KeTe04-SICON}, and \cite{Riff02-SICON}.

% -----------------------------------------------------------------------------------------------------------------------------------------
%		INSTABILITY THEOREMS
% -----------------------------------------------------------------------------------------------------------------------------------------

\section{Instability theorems}
\label{sec:Instability}
In addition to the stability considerations so far described, Lyapunov proposed using similar
energy-inspired functions for the study of instability.   In intuitive terms, if the system
energy is described by a positive function with a minimum at the origin, then stability, or asymptotic stability,
follows from the system energy not increasing, or decreasing, respectively.  This is clearly captured by the Lyapunov
functions discussed in previous sections.  In the case of instability, one sufficient condition has the
system energy (as described by a Lyapunov function) increasing in a neighborhood of the equilibrium
point.  However, this can be refined to allow that the system energy is increasing for points arbitrarily close
to the equilibrium point, but not necessarily in an entire neighborhood.  This is the underlying premise
of Chetaev's refinement to Lyapunov's instability theorems.
We first present the three instability theorems and then describe their converses.  The material
in this section is drawn from \cite{Kras63} and \cite{Hahn67}.

An unstable equilibrium is one that is not stable. {\rev
\begin{definition}
	The origin is unstable for \eqref{eq:sys} if, for any sufficiently small $\varepsilon > 0$ there exist sequences $\{x_k\}_{k=0}^\infty$ and
	$\{t_k\}_{k=0}^\infty$ such that $x_k \in \R^n \backslash \{0\}$ for all $k$, $x_k \rightarrow 0$ as $k \rightarrow \infty$,
	$t_k > t_0$ for all $k$, and
	\[ |\phi(t_k, t_0, x_k)| \geq \varepsilon, \quad \forall k . \]
\end{definition} }

Note that an equilibrium point can be both unstable and attractive; i.e., systems exist such that solutions from initial
conditions arbitrarily close to the origin leave every small neighborhood of the origin but eventually approach
the origin (see \cite[Section 40]{Hahn67}).

\begin{definition}
	Let $\G \subset \R^n$ contain the origin.  We say that the origin is {\em unstable in the region $\G$
	for \eqref{eq:sys}} if, for every open, bounded set $H \subset \R^n$ satisfying $0 \in H$ and $\overline{H} \subset \G$,
	and for every $t_0 \in \R_{\geq 0}$
	there exists a sequence of points $\{x_k\}_{k=0}^\infty$ satisfying $x_k \in H$, $\lim_{k \rightarrow \infty} x_k = 0$,
	and $\phi(t,t_0,x_k) \notin H$ for some $t > t_0$.
\end{definition}

In Lyapunov's First Theorem on Instability
\cite[Section 16, Theorem II]{Lyap1892} Lyapunov demonstrated that if a sign-definite function is such that, in a neighborhood
containing the origin,
its time derivative along solutions of \eqref{eq:sys} is also sign-definite and of the same sign
as the function itself, then the origin is unstable.
\begin{theorem}[Lyapunov's First Theorem on Instability]
	\label{thm:LyapunovFirstInstability}
	Let $\G \subset \R^n$ contain a neighborhood of the origin.
	Suppose the function $V: \G \times \R_{\geq 0} \rightarrow \R$ is continuously differentiable, is
	such that the derivative of $V$ along solutions of \eqref{eq:sys} is positive definite, and that there exists
	$\alpha \in \K$ such that $V(x,t) \leq \alpha(|x|)$ for all $(x,t) \in \G \times \R_{\geq 0}$.  If for every
	$\varepsilon > 0$ and $t_0 \geq 0$ there exists a $T \geq t_0$ such that, for all $|x| \leq \varepsilon$, $x \neq 0$,
	$V(x,t) > 0$ for all $t \geq T$, then the origin is unstable.
\end{theorem}

It is possible to obtain a converse to Theorem \ref{thm:LyapunovFirstInstability} if the origin is unstable
in the region $\G$ and $\G$ satisfies Property A (see \cite[Theorem 6.1]{Kras63}).

 In Lyapunov's second theorem on instability, \cite[Section 16, Theorem III]{Lyap1892},
 an extra degree of freedom is allowed in the form of a second function in the decrease condition.
\begin{theorem}[Lyapunov's Second Theorem on Instability]
\label{thm:LyapunovSecondInstability}
	Let $\G \subset \R^n$ contain a neighborhood of the origin.
	Suppose $W: \G \times \R_{\geq 0} \rightarrow \R_{\geq 0}$.  Suppose
	$V: \G \times \R_{\geq 0} \rightarrow \R_{\geq 0}$ is continuously differentiable, that there
	exists an $L > 0$ such that $|V(x,t)| \leq L$ for all $(x,t) \in \G \times \R_{\geq 0}$, and the
	derivative of $V$ along solutions of \eqref{eq:sys} satisfies
	\begin{equation}
		\frac{dV}{dt} = \lambda V + W
	\end{equation}
	for some $\lambda > 0$.  Furthermore, if $W(x,t) = 0$ for all $(x,t) \in \G \times \R_{\geq 0}$
	assume that for every $\varepsilon > 0$ and $t_0 \geq 0$ there exists a $T \geq t_0$
	such that, for all $|x| \leq \varepsilon$,
	$V(x,t) > 0$ for all $t \geq T$.  Then the origin is unstable.
\end{theorem}

It is possible to derive a converse to Theorem \ref{thm:LyapunovSecondInstability}
without the requirement that the region $\G$ satisfy Property A (see
\cite[Theorem 7.2]{Kras63}).

Finally, Chetaev \cite[Theorem, p.\ 27]{Chet61} revised the above to only require the function $V$
to be sign definite in a region containing the origin rather than in a neighborhood containing the origin.
\begin{theorem}[Chetaev's Theorem]
	\label{thm:Chetaev}
	Suppose $V: \R^n \times \R_{\geq 0} \rightarrow \R$ is continuously differentiable.
	If there exists $\varepsilon > 0$ such that, for all
	$|x| \leq \varepsilon$, the derivative of $V$ along solutions of \eqref{eq:sys} satisfies
	$dV/dt > 0$ on the region where $V(x,t) > 0$ then the origin
	is unstable.
\end{theorem}

The first converse theorems for the above instability theorems were derived by Krasovskii in the case of autonomous
systems \cite{Kras54-PMM}.  General converses for the instability theorems were derived
independently\footnote{Vrko\v{c}'s manuscript \cite{Vrko55-CMJ} was submitted on 7 January 1955 while
Krasovskii's manuscript \cite{Kras56-PMM} was submitted on 3 May 1955.}
by Vrko\v{c} \cite{Vrko55-CMJ} and Krasovskii \cite{Kras56-PMM}.

It is possible to derive a result that is slightly stronger than the direct converse to Theorem \ref{thm:Chetaev}.
To do this, we require the following definition.
\begin{definition}
	Let $H \subset \R^n$ be a bounded sub-domain of $\G$ that
	contains the origin and satisfies $\overline{H} \subset \G$.  A set $I(t_0) \subset H$,
	depending on the initial time $t_0$, is called the {\em domain of instability in $H$ for $t = t_0$} if for
	each $x \in I(t_0)$ there exists a finite time $t^* \in [t_0, \infty)$ so that
	$\phi(t^*,t_0,x) \notin H$.
\end{definition}
Note that the domain of instability is an open set in $\R^n$.

The following converse of Theorem \ref{thm:Chetaev} includes the result that the region where
the function $V$ is positive coincides with the region of instability.
\begin{theorem}{\cite[Theorem 7.1]{Kras63}}
	\label{thm:ChetaevConverse}
	Let the origin be unstable in the region $\G$ and let $H \subset \R^n$ be a bounded region
	satisfying $0 \in \overline{H} \subset \G$.  Then there exists a function
	$V:\G \times \R_{\geq 0} \rightarrow \R$ such that
	\begin{enumerate}
		\item If $x \in \overline{H}$ is in the region $V(x,t) > 0$ for all $t \geq t_0$
			then the function $dV/dt$ is positive definite;
		\item The function $V$ is bounded and continuous in the region
			$\overline{H}$ and the partial derivatives $\partial V / \partial t$ and
			$\partial V / \partial x_i$ are bounded uniformly in time; and
		\item For every value of $t = t_0$ the region of instability $I(t_0)$ coincides with
			the region $V(x,t) > 0$.
	\end{enumerate}
\end{theorem}

As in the results of Kurzweil and Massera on uniform asymptotic stability, if \eqref{eq:sys} is periodic in $t$ or independent of $t$,
then the function of Theorem \ref{thm:ChetaevConverse} can also be chosen to be periodic in $t$
or independent of $t$, respectively (see \cite[p.\ 43]{Kras63}).

We note that
\cite[Theorem 79]{Zubo64} extends Lyapunov's second theorem on instability to dynamical systems
defined on metric spaces.

% -----------------------------------------------------------------------------------------------------------------------------------------
%		 CONCLUSIONS
% -----------------------------------------------------------------------------------------------------------------------------------------

\section{Concluding remarks}
\label{sec:Conclusions}
As we have seen, the converse question for Lyapunov's second method has been successfully answered
in a wide variety of contexts.  The study of the converse question has helped to clarify not only the relationship
between different stability concepts, but has helped to identify useful stability concepts as in the case
of the important role played by uniformity in various stability definitions.  The answers have proved important
in the study of robustness to various system perturbations such as persistent disturbances, time delays, and
in the role that sampled-data controllers can play in providing robust feedback.

Of necessity, we have restricted our attention to
certain specific topics, specifically converse theorems for (uniform asymptotic) stability of
differential and difference inclusions, where the results for differential and difference equations
can be recovered as special cases.
However, a contributing factor in the success and popularity of Lyapunov's second method has been its
applicability in many different contexts and, in many of those contexts, converse results are also available.
For example, Lyapunov's second method can be extended to the notion of Lagrange stability
\cite{BaRo98-MCSS}.  An approach to almost global asymptotic stability, refered to as a ``dual'' to
Lyapunov's second method, was presented in \cite{Rant01-SCL} with a converse theorem in
\cite{Rant02-CDC}.  So-called complete Lyapunov functions for dynamical systems with
multiple equilibria have been defined in
\cite{Conl78} with converse theorems provided in \cite[Chapter 2, Section 6.4]{Conl78} and \cite{Patr11-FEJDS}.
Lyapunov's second method has also been adapted to the study of stochastic systems with converse
theorems available for
stochastic differential equations \cite{Kush67-SICON}, \cite[Theorem 5.4, p. 153]{Khas12},
discrete-time multivalued systems \cite{SuTe13-Aut}, \cite{THS14-TAC}, and
random dynamical systems for asymptotic stability of random compact sets \cite{ArSc01-JDE}.

%For acknowledgements section, please don't number the section, please begin it with \section*{Acknowledgements}
\section*{Acknowledgments}
{\rev Several colleagues were kind enough to read various drafts of this manuscript.  In particular, I would like
to thank Timm Faulwasser, Rafal Goebel, Lars Gr\"une, Petar Kokotovi\'c, Bj\"orn R\"uffer, and Dragoslav \v{S}iljak,
as well as the anonymous reviewers, for helpful comments.  Colleagues at the University of Bayreuth and
the Australian National University were also kind enough to invite me to present seminars on the material
contained herein and their feedback is also gratefully acknowledged.}

% You may incorporate your references as follows in your main tex file.
% Using BibTex is not recommended but can be handled.

\medskip
% The data information below will be filled by AIMS editorial staff
Received   August 2014; revised March 2015.
\medskip


\begin{thebibliography}{100}

\bibitem{Ande66-JFI} %(MR0204173) [10.1016/0016-0032(66)90317-6]
\newblock B.~D.~O. Anderson,
\newblock {Stability of control systems with multiple nonlinearities},
\newblock \emph{Journal of the Franklin Institute}, \textbf{282} (1966), 155--160.

\bibitem{AnMo69-SJC} %(MR0265020) [10.1137/0307029]
\newblock B.~D.~O. Anderson and J.~B. Moore,
\newblock {New results in linear system stability},
\newblock \emph{SIAM Journal on Control}, \textbf{7} (1969), 398--414.

\bibitem{AnMo81-SICON} %(MR603077) [10.1137/0319002]
\newblock B.~D.~O. Anderson and J.~B. Moore,
\newblock {Detectability and stabilizability of time-varying discrete-time linear systems},
\newblock \emph{{SIAM} Journal on Control and Optimization}, \textbf{19} (1981), 20--32.

\bibitem{AnSo99-SCL} %(MR1754903) [10.1016/S0167-6911(99)00055-9]
\newblock D.~Angeli and E.~D. Sontag,
\newblock {Forward completeness, unboundedness observability, and their {L}yapunov characterizations},
\newblock \emph{Systems \& Control Letters}, \textbf{38} (1999), 209--217.

\bibitem{Anto58} %(MR0102643)
\newblock H.~Antosiewicz,
\newblock A survey of {L}yapunov's second method,
\newblock \emph{Contributions to Nonlinear Oscillations}, Princeton University Press, 1958, 147--166.

\bibitem{Apos57} %(MR0087718)
\newblock T.~M. Apostol,
\newblock \emph{Mathematical Analysis: {A} Modern Approach to Advanced Calculus},
\newblock Addison-Wesley Publishing Company, Inc., Reading, Massachusetts, USA, 1957.

\bibitem{ArSc01-JDE} %(MR1867618) %[10.1006/jdeq.2000.3991]
\newblock L.~Arnold and B.~Schmalfuss,
\newblock {Lyapunov's second method for random dynamical systems},
\newblock \emph{Journal of Differential Equations}, \textbf{177} (2001), 235--265.

\bibitem{BaRo98-MCSS} %(MR1628047) %[10.1007/BF02741887]
\newblock A.~Bacciotti and L.~Rosier,
\newblock {Liapunov and {L}agrange stability: {I}nverse theorems for discontinuous systems},
\newblock \emph{Mathematics of Control, Signals and Systems}, \textbf{11} (1998), 101--128.

\bibitem{Barb48-Russian}
\newblock E.~A. Barbashin,
\newblock On the theory of general dynamical systems,
\newblock (Russian) {Ucen. Zap. Moskov. Gos. Univ.}, \textbf{135} (1948), 110--133.

\bibitem{Barb50-DAN} %(MR0036911)
\newblock E.~A. Barbashin,
\newblock Existence of smooth solutions of some linear equations with partial derivatives,
\newblock \emph{Doklady Akademii Nauk SSSR}, \textbf{72} (1950), 445--447.

\bibitem{BaKr52-DAN} %(MR0052616)
\newblock E.~A. Barbashin and N.~N. Krasovskii,
\newblock On the stability of motion in the large,
\newblock (Russian) \emph{Doklady Akademii Nauk SSSR}, \textbf{86} (1952), 453--456.

\bibitem{BaKr54-PMM} %(MR0062301)
\newblock E.~A. Barbashin and N.~N. Krasovskii,
\newblock On the existence of a function of {L}yapunov in the case of asymptotic stability in the large,
\newblock (Russian) \emph{Prikladnaya Matematika i Mekhanika}, \textbf{18} (1954), 345--350.

\bibitem{Broc82-BCh} %(MR708502)
\newblock R.~W. Brockett,
\newblock Asymptotic stability and feedback stabilization,
\newblock in \emph{Differential Geometric Control Theory} (eds. R.~W. Brockett, R.~S. Millman and H.~J. Sussman), Progr. Math., 27, Birkh\"auser Boston, Boston, MA, 1983, 181--191.

\bibitem{CTG07-TAC} %(MR2332751) [10.1109/TAC.2007.900829]
\newblock C.~Cai, A.~R. Teel and R.~Goebel,
\newblock {Smooth {L}yapunov functions for hybrid systems, {P}art {I}: {E}xistence is equivalent to robustness},
\newblock \emph{IEEE Transactions on Automatic Control}, \textbf{52} (2007), 1264--1277.

\bibitem{CTG08-TAC} %(MR2401025) [10.1109/TAC.2008.919257]
\newblock C.~Cai, A.~R. Teel and R.~Goebel,
\newblock {Smooth {L}yapunov functions for hybrid systems, {P}art {II}: {(Pre-)}asymptotically stable compact sets},
\newblock \emph{IEEE Transactions on Automatic Control}, \textbf{53} (2007), 734--748.

\bibitem{Chet61} %(MR0271487)
\newblock N.~G. Chetayev,
\newblock \emph{The Stability of Motion},
\newblock Pergamon Press, 1961; Translated from the 2nd Edition in Russian of 1956.

\bibitem{CLRS00-SICON} %(MR1780907) [10.1137/S0363012999352297]
\newblock F.~H. Clarke, Y.~S. Ledyaev, L.~Rifford and R.~J. Stern,
\newblock {Feedback stabilization and {L}yapunov functions},
\newblock \emph{SIAM Journal on Control and Optimization}, \textbf{39} (2000), 25--48.

\bibitem{CLSS97-TAC} %(MR1472857) [10.1109/9.633828]
\newblock F.~H. Clarke, Y.~S. Ledyaev, E.~D. Sontag and A.~I. Subbotin,
\newblock {Asymptotic controllability implies feedback stabilization},
\newblock \emph{IEEE Transactions on Automatic Control}, \textbf{42} (1997), 1394--1407.

\bibitem{CLS98-JDE} %(MR1643670) [10.1006/jdeq.1998.3476]
\newblock F.~H. Clarke, Y.~S. Ledyaev and R.~J. Stern,
\newblock {Asymptotic stability and smooth {L}yapunov functions},
\newblock \emph{Journal of Differential Equations}, \textbf{149} (1998), 69--114.

\bibitem{CLSW98} %(MR1488695)
\newblock F.~H. Clarke, Y.~S. Ledyaev, R.~J. Stern and P.~R. Wolenski,
\newblock \emph{Nonsmooth Analysis and Control Theory},
\newblock Springer-Verlag, New York, 1998.

\bibitem{Conl78} %(MR511133)
\newblock C.~Conley,
\newblock \emph{Isolated Invariant Sets and the {M}orse Index},
\newblock CBMS Regional Conference Series no.~38, American Mathematical Society, 1978.

\bibitem{CoTh06} %(MR2239987)
\newblock T.~M. Cover and J.~A. Thomas,
\newblock \emph{Elements of Information Theory},
\newblock 2nd edition, Wiley-Interscience, 2006.

\bibitem{Deim92} %(MR1189795) [10.1515/9783110874228]
\newblock K.~Deimling,
\newblock \emph{Multivalued Differential Equations},
\newblock Walter de Gruyter, 1992.

\bibitem{Fili88} %(MR1028776) [10.1007/978-94-015-7793-9]
\newblock A.~F. Filippov,
\newblock \emph{Differential Equations with Discontinuous Righthand Sides},
\newblock Kluwer Academic Publishers, 1988.

\bibitem{GiHa15-DCDSB}
\newblock P.~Giesl and S.~Hafstein,
\newblock Review on computational methods for {L}yapunov functions,
\newblock \emph{Discrete and Continuous Dynamical Systems, Series B}, \textbf{20} (2015).

\bibitem{GST12} %(MR2918932)
\newblock R.~Goebel, R.~G. Sanfelice and A.~R. Teel,
\newblock \emph{Hybrid Dynamical Systems: Modeling, Stability, and Robustness},
\newblock Princeton University Press, 2012.

\bibitem{Gord72-SJC} %(MR0318707) [10.1137/0310007]
\newblock S.~P. Gordon,
\newblock {On converses to the stability theorems for difference equations},
\newblock \emph{SIAM Journal on Control}, \textbf{10} (1972), 76--81.

\bibitem{GCW01-SICON} %(MR1857360) [10.1137/S036301299936316X]
\newblock L.~Gr\"une, F.~Camilli and F.~Wirth,
\newblock {A generalization of {Z}ubov's method to perturbed systems},
\newblock \emph{SIAM Journal on Control and Optimization}, \textbf{40} (2001), 496--515.

\bibitem{GKSW07-DCDS} %(MR2291904) [10.3934/dcds.2007.18.375]
\newblock L.~Gr\"une, P.~E. Kloeden, S.~Siegmund and F.~R. Wirth,
\newblock {Lyapunov's second method for nonautonomous differential equations},
\newblock \emph{Discrete and Continuous Dynamical Systems}, \textbf{18} (2007), 375--403.

\bibitem{GrSe11-SICON} %(MR2854621) [10.1137/100787829]
\newblock L.~Gr\"une and O.~S. Serea,
\newblock {Differential games and {Z}ubov's method},
\newblock \emph{SIAM Journal on Control and Optimization}, \textbf{49} (2011), 2349--2377.

\bibitem{Hahn63} %(MR0147716)
\newblock W.~Hahn,
\newblock \emph{Theory and Application of Liapunov's Direct Method},
\newblock Prentice-Hall, 1963; Translated from the German Edition of 1959.

\bibitem{Hahn67} %(MR0223668)
\newblock W.~Hahn,
\newblock \emph{Stability of Motion},
\newblock Springer-Verlag, 1967.

\bibitem{HiAn69-IEE} %(MR0309612) [10.1049/piee.1969.0031]
\newblock B.~E. Hitz and B.~D.~O. Anderson,
\newblock {Discrete positive-real functions and their application to system stability},
\newblock \emph{Proceedings of the Institution of Electrical Engineers}, \textbf{116} (1969), 153--155.

\bibitem{Hopp66-TAMS} %(MR0194693) [10.1090/S0002-9947-1966-0194693-9]
\newblock F.~C. Hoppensteadt,
\newblock {Singular perturbations on the infinite interval},
\newblock \emph{Transactions of the American Mathematical Society}, \textbf{123} (1966), 521--535.

\bibitem{ISW02-CDC} %[10.1109/CDC.2002.1184983]
\newblock B.~P. Ingalls, E.~D. Sontag and Y.~Wang,
\newblock {Measurement to error stability: A notion of partial detectability for nonlinear systems},
\newblock in \emph{Proceedings of the 41st IEEE Conference on Decision and Control, Vol. 4}, Las Vegas, Nevada, USA, 2002, 3946--3951.

\bibitem{JiWa02-SCL} %(MR2010491) [10.1016/S0167-6911(01)00164-5]
\newblock Z.-P. Jiang and Y.~Wang,
\newblock {A converse {L}yapunov theorem for discrete-time systems with disturbances},
\newblock \emph{Systems \& Control Letters}, \textbf{45} (2002), 49--58.

\bibitem{Kalm63-PNAS} %(MR0151696) [10.1073/pnas.49.2.201]
\newblock R.~E. Kalman,
\newblock {Lyapunov function for the problem of {L}ur'e in automatic control},
\newblock \emph{Proc. Nat. Acad. Sci. U.S.A.}, \textbf{49} (1963), 201--205.

\bibitem{KaBe60-JBE} %(MR0157810) [10.1115/1.3662604]
\newblock R.~E. Kalman and J.~E. Bertram,
\newblock {Control system analysis and design via the ``second method'' of {L}yapunov, {P}art {I}, continuous-time systems},
\newblock \emph{Transactions of the AMSE, Series D: Journal of Basic Engineering}, \textbf{82} (1960), 371--393.

\bibitem{KaBe60b-JBE} %(MR0157811) [10.1115/1.3662605]
\newblock R.~E. Kalman and J.~E. Bertram,
\newblock {Control system analysis and design via the ``second method'' of {L}yapunov, {P}art {II}, discrete-time systems},
\newblock \emph{Transactions of the AMSE, Series D: Journal of Basic Engineering}, \textbf{82} (1960), 394--400.

\bibitem{Kara05-IMA} %(MR2212259) [10.1093/imamci/dni037]
\newblock I.~Karafyllis,
\newblock {Non-uniform robust global asymptotic stability for discrete-time systems and applications to numerical analysis},
\newblock \emph{IMA Journal of Mathematical Control and Information}, \textbf{23} (2006), 11--41.

\bibitem{KaTs03-SICON} %(MR2002141) [10.1137/S0363012901392967]
\newblock I.~Karafyllis and J.~Tsinias,
\newblock {A converse {L}yapunov theorem for nonuniform in time global asymptotic stability and its application to feedback stabilization},
\newblock \emph{{SIAM} Journal on Control and Optimization}, \textbf{42} (2003), 936--965.

\bibitem{Kell14-MCSS} %(MR3245919) [10.1007/s00498-014-0128-8]
\newblock C.~M. Kellett,
\newblock {A compendium of comparison function results},
\newblock \emph{Mathematics of Controls, Signals and Systems}, \textbf{26} (2014), 339--374.

\bibitem{KeTe00-MTNS}
\newblock C.~M. Kellett and A.~R. Teel,
\newblock A converse {L}yapunov theorem for weak uniform asymptotic stability of sets,
\newblock in \emph{Proceedings of Mathematical Theory of Networks and Systems}, Perpignan, France, 2000.

\bibitem{KeTe00-CDC} %[10.1109/CDC.2000.912339]
\newblock C.~M. Kellett and A.~R. Teel,
\newblock {Uniform asymptotic controllability to a set implies locally {L}ipschitz control-{L}yapunov function},
\newblock in \emph{Proceedings of the 39th IEEE Conference on Decision and Control, Vol. 4}, Sydney, Australia, 2000, 3994--3999.

\bibitem{KeTe04b-SCL} %(MR2074378) [10.1016/j.sysconle.2004.02.011]
\newblock C.~M. Kellett and A.~R. Teel,
\newblock {Discrete-time asymptotic controllability implies smooth control-{L}yapunov function},
\newblock \emph{Systems \& Control Letters}, \textbf{52} (2004), 349--359.

\bibitem{KeTe04-SICON} %(MR2080922) [10.1137/S0363012901398186]
\newblock C.~M. Kellett and A.~R. Teel,
\newblock {Weak converse {L}yapunov theorems and control {L}yapunov functions},
\newblock \emph{SIAM Journal on Control and Optimization}, \textbf{42} (2004), 1934--1959.

\bibitem{KeTe05-SICON} %(MR2178046) [10.1137/S0363012903435862]
\newblock C.~M. Kellett and A.~R. Teel,
\newblock {On the robustness of $\mathcal{KL}$-stability for difference inclusions: {S}mooth discrete-time {L}yapunov functions},
\newblock \emph{SIAM Journal on Control and Optimization}, \textbf{44} (2005), 777--800.

\bibitem{KeTe07-MCSS} %(MR2324813) [10.1007/s00498-007-0016-6]
\newblock C.~M. Kellett and A.~R. Teel,
\newblock {Sufficient conditions for robustness of $\mathcal{KL}$-stability for difference inclusions},
\newblock \emph{Mathematics of Control, Signals and Systems}, \textbf{19} (2007), 183--205.

\bibitem{Khal96}
\newblock H.~K. Khalil,
\newblock \emph{Nonlinear Systems},
\newblock 2nd edition, Prentice Hall, 1996.

\bibitem{Khas12} %(MR2894052) [10.1007/978-3-642-23280-0]
\newblock R.~Khasminskii,
\newblock \emph{Stochastic Stability of Differential Equations},
\newblock 2nd edition, Springer, 2012.

\bibitem{Kloe78-MCT} %(MR515715)
\newblock P.~E. Kloeden,
\newblock General control systems,
\newblock in \emph{Mathematical Control Theory 1977: Proceedings} (ed. W.~A. Coppel), Springer-Verlag, Canberra, Australia, 1978, 119--137.

\bibitem{Kloe98} %(MR1659222)
\newblock P.~E. Kloeden,
\newblock Lyapunov functions for cocycle attractors in nonautonomous difference equations,
\newblock \emph{Izvetsiya Akad Nauk Rep Moldovia Mathematika}, \textbf{26} (1998), 32--42.

\bibitem{Kloe00} %(MR1799047)
\newblock P.~E. Kloeden,
\newblock A {L}yapunov function for pullback attractors of nonautonomous differential equations,
\newblock \emph{Electronic Journal of Differential Equations Conference 05}, (2000), 91--102.

\bibitem{KoAr01-Aut} %(MR1832954) [10.1016/S0005-1098(01)00002-4]
\newblock P.~Kokotovi\'{c} and M.~Arcak,
\newblock {Constructive nonlinear control: A historical perspective},
\newblock \emph{Automatica}, \textbf{37} (2001), 637--662.

\bibitem{Kras54-PMM} %(MR0065729)
\newblock N.~N. Krasovskii,
\newblock On the inversion of the theorems of {A.\ M.\ Lyapunov} and {N.\ G.\ Chetaev} concerning instability for stationary systems of differential equations,
\newblock (Russian) \emph{Prikladnaya Matematika i Mekhanika}, \textbf{18} (1954), 513--532.

\bibitem{Kras55-PMM} %(MR0069348)
\newblock N.~N. Krasovskii,
\newblock On the converse of {K. P. Persidskii's} theorem on uniform stability,
\newblock (Russian) \emph{Prikladnaya Matematika i Mekhanika}, \textbf{19} (1955), 273--278.

\bibitem{Kras56-PMM}
\newblock N.~N. Krasovskii,
\newblock Transformation of the theorem of {A. M. Lyapunov's} second method and questions of first-order stability of motion,
\newblock (Russian) \emph{Prikladnaya Matematika i Mekhanika}, \textbf{20} (1956), 255--265.

\bibitem{Kras63} %(MR0147744)
\newblock N.~N. Krasovskii,
\newblock \emph{Stability of Motion: {A}pplications of {L}yapunov's Second Method to Differential Systems and Equations with Delay},
\newblock Stanford University Press, 1963; Translated from the Russian Edition of 1959.

\bibitem{KSW01-SICON} %(MR1825868) [10.1137/S0363012999365352]
\newblock M.~Krichman, E.~D. Sontag and Y.~Wang,
\newblock {Input-output-to-state stability},
\newblock \emph{{SIAM} Journal on Control and Optimization}, \textbf{39} (2001), 1874--1928.

\bibitem{KKK95}
\newblock M.~Krsti\'{c}, I.~Kanellakopoulos and P.~Kokotovi\'{c},
\newblock \emph{Nonlinear and Adaptive Control Design},
\newblock John Wiley and Sons, Inc., 1995.

\bibitem{Kurz55-CMJ} %(MR0077747)
\newblock J.~Kurzweil,
\newblock Transformation of {L}yapunov's first theorem on stability of motion,
\newblock (Russian) \emph{Czechoslovak Mathematical Journal}, \textbf{5} (1955), 382--398.

\bibitem{Kurz56-AMST}
\newblock J.~Kurzweil,
\newblock On the inversion of {L}japunov's second theorem on stability of motion,
\newblock (Russian) \emph{Czechoslovak Mathematical Journal}, \textbf{81} (1956), 217--259, 455--484; English translation in \emph{American Mathematical Society Translations (2)}, \textbf{24}, 19--77.

\bibitem{KuVr57-CMJ} %(MR0089325)
\newblock J.~Kurzweil and I.~Vrko\v{c},
\newblock Transformation of {L}yapunov's theorems on stability and {P}ersidskii's theorems on uniform stability,
\newblock (Russian) \emph{Czechoslovak Mathematical Journal}, \textbf{7} (1957), 254--272.

\bibitem{Kush67-SICON} %(MR0213183) [10.1137/0305015]
\newblock H.~J. Kushner,
\newblock {Converse theorems for stochastic {L}iapunov functions},
\newblock \emph{SIAM Journal on Control and Optimization}, \textbf{5} (1967), 228--233.

\bibitem{LaLe61} %(MR0132876)
\newblock J.~{La Salle} and S.~Lefschetz,
\newblock \emph{Stability by {L}iapunov's Direct Method with Applications},
\newblock Academic Press, 1961.

\bibitem{LaSa76-BUMI} %(MR0440136)
\newblock V.~Lakshmikantham and L.~Salvadori,
\newblock On {M}assera type converse theorem in terms of two different measures,
\newblock \emph{Bollettino dell'Unione Matematica Italiana}, \textbf{13} (1976), 293--301.

\bibitem{Leto61} %(MR0129081)
\newblock A.~L. Letov,
\newblock \emph{Stability in Nonlinear Control Systems},
\newblock Princeton University Press, Princeton, New Jersey, 1961; Translated from the Russian Edition of 1955.

\bibitem{LSW96-SICON} %(MR1372908) [10.1137/S0363012993259981]
\newblock Y.~Lin, E.~D. Sontag and Y.~Wang,
\newblock {A smooth converse {L}yapunov theorem for robust stability},
\newblock \emph{{SIAM} Journal on Control and Optimization}, \textbf{34} (1996), 124--160.

\bibitem{Lure51}
\newblock A.~I. Lur'e,
\newblock \emph{Some Non-Linear Problems in the Theory of Automatic Control},
\newblock Her Majesty's Stationery Office, 1957; Translated from the Russian Edition of 1951.

\bibitem{LuPo44-PMM} %(MR0011360)
\newblock A.~I. Lur'e and V.~N. Postnikov,
\newblock Stability theory of regulating systems,
\newblock (Russian) \emph{Prikladnaya Matematika i Mekhanika}, \textbf{8} (1944), 246--248.

\bibitem{Lyap1892} %(MR1154209) [10.1080/00207179208934253]
\newblock A.~M. Lyapunov,
\newblock {The general problem of the stability of motion},
\newblock (Russian) \emph{Math. Soc. of Kharkov}; English Translation, \emph{International Journal of Control}, \textbf{55} (1992), 531--773.

\bibitem{Malk54-PMM}
\newblock I.~G. Malkin,
\newblock Questions concerning transformation of {L}yapunov's theorem on asymptotic stability,
\newblock (Russian) \emph{Prikladnaya Matematika i Mekhanika}, \textbf{18} (1954), 129--138.

\bibitem{Malk56}
\newblock I.~G. Malkin,
\newblock \emph{Some Problems in the Theory of Nonlinear Oscillations},
\newblock United States Atomic Energy Commission, 1959, Translated from the Russian Edition of 1956.

\bibitem{Mass49-AM} %(MR0035354) [10.2307/1969558]
\newblock J.~L. Massera,
\newblock {On {L}iapounoff's conditions of stability},
\newblock \emph{Annals of Mathematics}, \textbf{50} (1949), 705--721.

\bibitem{Mass56-AM} %(MR0079179) [10.2307/1969955]
\newblock J.~L. Massera,
\newblock {Contributions to stability theory},
\newblock \emph{Annals of Mathematics}, \textbf{64} (1956), 182--206.
\newblock (Erratum: \textit{Annals of Mathematics}, \textbf{68} (1958), 202.)

\bibitem{Meil78-AiT} %(MR533365)
\newblock A.~M. Meilakhs,
\newblock Design of stable control systems subject to parametric perturbation,
\newblock \emph{Automation and Remote Control}, \textbf{39} (1979), 1409--1418.

\bibitem{MHL08} %(MR2351563)
\newblock A.~N. Michel, L.~Hou and D.~Liu,
\newblock \emph{Stability of Dynamical Systems: Continuous, Discontinuous, and Discrete Systems},
\newblock Birkh\"auser, 2008.

\bibitem{MoPy86a-AiT} %(MR839959)
\newblock A.~P. Molchanov and E.~S. Pyatnitskii,
\newblock Lyapunov functions that specifiy necessary and sufficient conditions of absolute stability of nonlinear nonstationary control systems {I},
\newblock  \emph{Automation and Remote Control}, \textbf{47} (1986), 344--354.

\bibitem{MoPy86b-AiT} %(MR848396)
\newblock A.~P. Molchanov and E.~S. Pyatnitskii,
\newblock Lyapunov functions that specifiy necessary and sufficient conditions of absolute stability of nonlinear nonstationary control systems {II},
\newblock  \emph{Automation and Remote Control}, \textbf{47} (1986), 443--451.

\bibitem{MoPy86c-AiT} %(MR848396)
\newblock A.~P. Molchanov and E.~S. Pyatnitskii,
\newblock Lyapunov functions that specifiy necessary and sufficient conditions of absolute stability of nonlinear nonstationary control systems {III},
\newblock \emph{Automation and Remote Control}, \textbf{47} (1986), 620--630.

\bibitem{MoPy89-SCL} %(MR1006848) [10.1016/0167-6911(89)90021-2]
\newblock A.~P. Molchanov and Y.~S. Pyatnitskiy,
\newblock {Criteria of asymptotic stability of differential and difference inclusions encountered in control theory},
\newblock \emph{Systems \& Control Letters}, \textbf{13} (1989), 59--64.

\bibitem{Movc60-PMM} %[10.1016/0021-8928(60)90004-6]
\newblock A.~A. Movchan,
\newblock {Stability of processes with respect to two metrics},
\newblock \emph{Journal of Applied Mathematics and Mechanics}, \textbf{24} (1960), 1506--1524.

\bibitem{Patr11-FEJDS} %(MR2934472)
\newblock M.~Patrao,
\newblock Existence of complete {L}yapunov functions for semiflows on separable metric spaces,
\newblock \emph{Far East Journal of Dynamical Systems}, \textbf{17} (2011), 49--54.

\bibitem{Pers37-Doklady}
\newblock K.~P. Persidskii,
\newblock On a theorem of {L}iapunov,
\newblock \emph{C. R. (Dokl.) Acad. Sci. URSS}, \textbf{14} (1937), 541--543.

\bibitem{Popo63-AiT} %(MR0133563)
\newblock V.~M. Popov,
\newblock Absolute stability of nonlinear systems of automatic control,
\newblock \emph{Automation and Remote Control}, \textbf{22} (1961), 857--875.

\bibitem{Popo64}
\newblock V.~M. Popov,
\newblock Proprietati de stabilitate si de optimalitate pentru sistemele
  automate cu mai multe functii de comanda,
\newblock (Romanian) \emph{Studii si Cercetari de Energetica, Academici RPR}, \textbf{14}
  (1964), 913--949.

\bibitem{Popo66} %(MR0387749)
\newblock V.~M. Popov,
\newblock \emph{Hyperstability of Control Systems},
\newblock Springer-Verlag, 1973, Translated from the Romanian Edition of 1966.

\bibitem{Rant01-SCL} %(MR2007046) [10.1016/S0167-6911(00)00087-6]
\newblock A.~Rantzer,
\newblock {A dual to {L}yapunov's stability theorem},
\newblock \emph{Systems \& Control Letters}, \textbf{42} (2001), 161--168.

\bibitem{Rant02-CDC} %[10.1109/CDC.2002.1184801]
\newblock A.~Rantzer,
\newblock {An converse theorem for density functions},
\newblock in \emph{Proceedings of the 41st IEEE Conference on Decision and Control, Vol. 2}, Las Vegas, Nevada, USA, 2002, 1890--1891.

\bibitem{Riff00-SICON} %(MR1814266) [10.1137/S0363012999356039]
\newblock L.~Rifford,
\newblock {Existence of {L}ipschitz and semiconcave control-{L}yapunov functions},
\newblock \emph{SIAM Journal on Control and Optimization}, \textbf{39} (2000), 1043--1064.

\bibitem{Riff02-SICON} %(MR1939865) [10.1137/S0363012900375342]
\newblock L.~Rifford,
\newblock {Semiconcave control-{L}yapunov functions and stabilizing feedbacks},
\newblock \emph{SIAM Journal on Control and Optimization}, \textbf{41} (2002), 659--681.

\bibitem{Rosi92-SCL} %(MR1195304) [10.1016/0167-6911(92)90078-7]
\newblock L.~Rosier,
\newblock {Homogeneous {L}yapunov function for homogeneous continuous vector field},
\newblock \emph{Systems \& Control Letters}, \textbf{19} (1992), 467--473.

\bibitem{RHL77} %(MR0450715)
\newblock N.~Rouche, P.~Habets and M.~Laloy,
\newblock \emph{Stability Theory by Liapunov's Direct Method},
\newblock Springer-Verlag, 1977.

\bibitem{Roxi65b-JDE} %(MR0201756) [10.1016/0022-0396(65)90019-7]
\newblock E.~Roxin,
\newblock {On generalized dynamical systems defined by contingent equations},
\newblock \emph{Journal of Differential Equations}, \textbf{1} (1965), 188--205.

\bibitem{Roxi65a-JDE} %(MR0201755) [10.1016/0022-0396(65)90015-X]
\newblock E.~Roxin,
\newblock {Stability in general control systems},
\newblock \emph{Journal of Differential Equations}, \textbf{1} (1965), 115--150.

\bibitem{Roxi66-RCMP} %(MR0230995) [10.1007/BF02849435]
\newblock E.~Roxin,
\newblock {On asymptotic stability in control systems},
\newblock \emph{Rendiconti del Circolo Matematico di Palermo}, \textbf{15} (1966), 193--208.

\bibitem{Roxi66-SJC} %(MR0196243) [10.1137/0303024]
\newblock E.~Roxin,
\newblock {On stability in control systems},
\newblock \emph{SIAM Journal on Control}, \textbf{3} (1966), 357--372.

\bibitem{ShWe49} %(MR0032134)
\newblock C.~E. Shannon and W.~Weaver,
\newblock \emph{The Mathematical Theory of Communication},
\newblock University of Illinois Press, 1949.

\bibitem{SWMWK07-SIAM-Review} %(MR2375524) [10.1137/05063516X]
\newblock R.~Shorten, F.~Wirth, O.~Mason, K.~Wulff and C.~King,
\newblock {Stability criteria for switched and hybrid systems},
\newblock \emph{SIAM Review}, \textbf{49} (2007), 545--592.

\bibitem{Silj69}
\newblock D.~D. \v{S}iljak,
\newblock \emph{Nonlinear Systems: The Parameter Analysis and Design},
\newblock John Wiley \& Sons Inc., 1969.

\bibitem{Smir90a-AiT} %(MR1071217)
\newblock G.~V. Smirnov,
\newblock Weak asymptotic stability of differential inclusions {I},
\newblock \emph{Automation and Remote Control}, \textbf{51} (1990), 901--908.

\bibitem{Smir90b-AiT} %(MR1080602)
\newblock G.~V. Smirnov,
\newblock Weak asymptotic stability of differential inclusions {II},
\newblock \emph{Automation and Remote Control}, \textbf{51} (1990), 1052--1058.

\bibitem{Sont83-SICON} %(MR696908) [10.1137/0321028]
\newblock E.~D. Sontag,
\newblock {A {L}yapunov-like characterization of asymptotic controllability},
\newblock \emph{{SIAM} Journal on Control and Optimization}, \textbf{21} (1983), 462--471.

\bibitem{Sont89-TAC} %(MR987806) [10.1109/9.28018]
\newblock E.~D. Sontag,
\newblock {Smooth stabilization implies coprime factorization},
\newblock \emph{IEEE Transactions on Automatic Control}, \textbf{34} (1989), 435--443.

\bibitem{Sont99-COCV} %(MR1746166) [10.1051/cocv:1999121]
\newblock E.~D. Sontag,
\newblock {Clocks and insensitivity to small measurement errors},
\newblock \emph{{ESAIM}: Control, Optimization, and the Calculus of Variations}, \textbf{4} (1999), 537--557.

\bibitem{SoWa95-SCL} %(MR1325675) [10.1016/0167-6911(94)00050-6]
\newblock E.~D. Sontag and Y.~Wang,
\newblock {On characterizations of the input-to-state stability property},
\newblock \emph{Systems \& Control Letters}, \textbf{24} (1995), 351--359.

\bibitem{Stei52-JRNBS} %(MR0047001) [10.6028/jres.048.010]
\newblock P.~Stein,
\newblock {Some general theorems on iterants},
\newblock \emph{Journal of Research of the National Bureau of Standards}, \textbf{48} (1952), 82--83.

\bibitem{SuTe13-Aut} %(MR3092645) [10.1016/j.automatica.2013.07.001]
\newblock A.~Subbaraman and A.~R. Teel,
\newblock {A converse {L}yapunov theorem for strong global recurrence},
\newblock \emph{Automatica}, \textbf{49} (2013), 2963--2974.

\bibitem{THS14-TAC} %(MR3254536) [10.1109/TAC.2014.2322431]
\newblock A.~R. Teel, J.~P. Hespanha and A.~Subbaraman,
\newblock {A converse {L}yapunov theorem and robustness for asymptotic stability in probability},
\newblock \emph{IEEE Transactions on Automatic Control}, \textbf{59} (2014), 2426--2421.

\bibitem{TePr00-CoCV} %(MR1765429) [10.1051/cocv:2000113]
\newblock A.~R. Teel and L.~Praly,
\newblock {A smooth {L}yapunov function from a class-$\mathcal{KL}$ estimate involving two positive semidefinite functions},
\newblock \emph{{ESAIM}: Control, Optimization, and the Calculus of Variations}, \textbf{5} (2000), 313--367.

\bibitem{Tsyp62-Doklady}
\newblock Y.~Z. Tsypkin,
\newblock The absolute stability of large-scale nonlinear sampled-data systems,
\newblock (Russian) \emph{Doklady Akademii Nauk SSSR}, \textbf{145} (1962), 52--55.

\bibitem{Tsyp63-AiT} %(MR0159712)
\newblock Y.~Z. Tsypkin,
\newblock Absolute stability of equilibrium positions and of responses in nonlinear, sampled-data automatic systems,
\newblock  \emph{Automation and Remote Control}, \textbf{24} (1963), 1457--1470.

\bibitem{Voro05-AiT} %(MR2140056) [10.1007/s10513-005-0099-9]
\newblock V.~I. Vorotnikov,
\newblock {Partial stability and control: {T}he state-of-the-art and development prospects},
\newblock  \emph{Automation and Remote Control}, \textbf{66} (2005), 511--561.

\bibitem{Vrko55-CMJ}
\newblock I.~Vrko\v{c},
\newblock A general theorem of {C}hetaev,
\newblock (Russian) \emph{Czechoslovak Mathematical Journal}, \textbf{5} (1955), 451--461.

\bibitem{Will-ARMA72a} %(MR0527462) [10.1007/BF00276493]
\newblock J.~C. Willems,
\newblock {Dissipative dynamical systems part {I}: {G}eneral theory},
\newblock \emph{Archive for Rational Mechanics and Analysis}, \textbf{45} (1972), 321--351.

\bibitem{Will-ARMA72b} %(MR0527463) [10.1007/BF00276494]
\newblock J.~C. Willems,
\newblock {Dissipative dynamical systems part {II}: {L}inear systems with quadratic supply rates},
\newblock \emph{Archive for Rational Mechanics and Analysis}, \textbf{45} (1972), 352--393.

\bibitem{Wils69-TAMS} %(MR0251747) [10.1090/S0002-9947-1969-0251747-9]
\newblock F.~W. Wilson,
\newblock {Smoothing derivatives of functions and applications},
\newblock \emph{Transactions of the American Mathematical Society}, \textbf{139} (1969), 413--428.

\bibitem{Yaku62-DAN} %(MR0133333)
\newblock V.~A. Yakubovich,
\newblock The solution of certain matrix inequalities in automatic control theory,
\newblock \emph{Doklady Akademii Nauk SSSR}, \textbf{143} (1962), 1304--1307.

\bibitem{Yosh55-Kyoto} %(MR0075383)
\newblock T.~Yoshizawa,
\newblock On the stability of solutions of a system of differential equations,
\newblock \emph{Mem. Coll. Sci. Univ. Kyoto. Ser. A. Math.}, \textbf{29} (1955), 27--33.

\bibitem{Yosh66} %(MR0208086)
\newblock T.~Yoshizawa,
\newblock \emph{Stability Theory by Liapunov's Second Method},
\newblock Mathematical Society of Japan, 1966.

\bibitem{Zubo64} %(MR0179428)
\newblock V.~I. Zubov,
\newblock \emph{Methods of A.~M.~Lyapunov and their Application},
\newblock P.~Noordhoff Ltd, Groningen, The Netherlands, 1964; Translated from the Russian Edition of 1957.

\end{thebibliography}
\end{document}